\documentclass[conference]{IEEEtran}
\IEEEoverridecommandlockouts
\usepackage{cite}
\usepackage{amsmath,amssymb,amsfonts}
\usepackage{algorithmic}
\usepackage{graphicx}
\usepackage{textcomp}
\usepackage{xcolor}
\def\BibTeX{{\rm B\kern-.05em{\sc i\kern-.025em b}\kern-.08em
    T\kern-.1667em\lower.7ex\hbox{E}\kern-.125emX}}
    
\usepackage[utf8]{inputenc}
\usepackage{pgfplots}
\DeclareUnicodeCharacter{2212}{−}
\usepgfplotslibrary{groupplots,dateplot}
\usetikzlibrary{patterns,shapes.arrows}
\pgfplotsset{compat=newest}
    
\newcommand{\bs}[1]{{\ensuremath{\boldsymbol{#1}}}}
\DeclareMathOperator*{\argmin}{\arg\!\min}

\usepackage{siunitx}
\usepackage[nolist]{acronym}

\usetikzlibrary{shapes.geometric}
\usetikzlibrary{arrows.meta,arrows}

\def\BibTeX{{\rm B\kern-.05em{\sc i\kern-.025em b}\kern-.08em
    T\kern-.1667em\lower.7ex\hbox{E}\kern-.125emX}}
    
\newcommand{\greencheck}{}%
\DeclareRobustCommand{\greencheck}{%
  \tikz\fill[scale=0.4, color=black!60!green]
  (0,.35) -- (.25,0) -- (1,.7) -- (.25,.15) -- cycle;%
}

\newcommand{\redcross}{}%
\DeclareRobustCommand{\redcross}{%
  \tikz[scale=0.23, color=red] {
    \draw[line width=0.7,line cap=round] (0,0) to [bend left=6] (1,1);
    \draw[line width=0.7,line cap=round] (0.2,0.95) to [bend right=3] (0.8,0.05);
}
}
\newcolumntype{R}[2]{%
    >{\adjustbox{angle=#1,lap=\width-(#2)}\bgroup}%
    l%
    <{\egroup}%
}

\renewcommand{\footnoterule}{%
  \kern -3pt
  \hrule width 6cm height 1pt
  \kern 2pt
}

\usepackage{eurosym}
\DeclareSIUnit{\sieuro}{\mbox{\euro}}
\begin{document}
\begin{acronym}
    \acro{CEA}{controlled environment agriculture}
    \acro{GH}{greenhouse}
    \acro{OCP}{optimal control problem}
    \acro{PPFD}{photosynthetic photon flux density}
    \acro{OC}{optimal control}
\end{acronym}
\title{Comparison of Dynamic Tomato Growth Models for Optimal Control in Greenhouses\\
}
\author{
\IEEEauthorblockN{ Michael Fink$^*$}
\IEEEauthorblockA{\small \textit{Chair of Automatic Control Engineering} \\
\textit{Technical University of Munich}, Germany\\
michael.fink@tum.de}
\and 
\IEEEauthorblockN{Annalena Daniels$^*$ }
\IEEEauthorblockA{\small \textit{Chair of Automatic Control Engineering} \\
\textit{Technical University of Munich}, Germany\\
a.daniels@tum.de}
\and
\IEEEauthorblockN{Cheng Qian}
\IEEEauthorblockA{\small \textit{Chair of Automatic Control Engineering} \\
\textit{Technical University of Munich}, Germany\\
cheng.qian@tum.de}
\and

\IEEEauthorblockN{Víctor Martínez Velásquez}
\IEEEauthorblockA{\small \textit{Chair of Automatic Control Engineering} \\
\textit{Technical University of Munich}, Germany\\
victor.martinez@tum.de}
\and 
\IEEEauthorblockN{Sahil Salotra}
\IEEEauthorblockA{\small \textit{Chair of Automatic Control Engineering} \\
\textit{Technical University of Munich}, Germany\\
sahil.salotra@tum.de}
\and
\IEEEauthorblockN{Dirk Wollherr}
\IEEEauthorblockA{\small \textit{Chair of Automatic Control Engineering} \\
\textit{Technical University of Munich}, Germany\\
dw@tum.de}
}

\maketitle
\def\thefootnote{*}\footnotetext{These authors contributed equally to this work. \\
© 2023 IEEE. This work was accepted for IEEE AGRETA2023. Personal use of this material is permitted.  Permission from IEEE must be obtained for all other uses, in any current or future media, including reprinting/republishing this material for advertising or promotional purposes, creating new collective works, for resale or redistribution to servers or lists, or reuse of any copyrighted component of this work in other works.
}\def\thefootnote{\arabic{footnote}}

\begin{abstract}
As global demand for efficiency in agriculture rises, there is a growing interest in high-precision farming practices. Particularly greenhouses play a critical role in ensuring a year-round supply of fresh produce. In order to maximize efficiency and productivity while minimizing resource use, mathematical techniques such as optimal control have been employed. However, selecting appropriate models for optimal control requires domain expertise. This study aims to compare three established tomato models for their suitability in an optimal control framework. Results show that all three models have similar yield predictions and accuracy, but only two models are currently applicable for optimal control due to implementation limitations. The two remaining models each have advantages in terms of economic yield and computation times, but the differences in optimal control strategies suggest that they require more accurate parameter identification and calibration tailored to greenhouses.
\end{abstract}


\acresetall
\section{Introduction}
With a growing world population, the agricultural sector is challenged to increase food production while minimizing the negative impact on the environment and preserving biodiversity \cite{Searchinger.2019}. 
In addition, the availability of agricultural land is decreasing due to factors including climate change and geopolitical conflicts, highlighting the need for higher production density. To achieve higher yields and quality while reducing cost and environmental impact, there is a growing trend towards high-precision \acl{CEA} \cite{Shamshiri.2018}, of which \acp{GH} are a crucial element. \acp{GH} are partially enclosed systems that regulate environmental variables like temperature, humidity and CO\textsubscript{2} concentration.

Despite their benefits for providing food all-year-round, \acp{GH} still face challenges due to their high energy and cost requirements, making them less profitable than arable farming during cropping seasons. Nonetheless, they offer significant advantages for crops like tomatoes, which are difficult to transport and store for extended periods. Maximizing the potential yield of these crops while minimizing the economic and ecological costs can unlock the full potential of \ac{GH} tomato production.

The optimization of \ac{GH} management for enhanced crop productivity and reduced resource use and environmental impact requires the use of \ac{OC}. \ac{OC} is a mathematical technique that aims at finding the best control strategy for a given system. It involves determining the optimal values of control variables over a specified time horizon that maximize or minimize a performance measure, subject to constraints. In the context of \ac{GH} management, \ac{OC} can be used to optimize environmental factors such as temperature, humidity, and lighting to enhance crop productivity while minimizing resource use and environmental impact \cite{vanStraten.2010, ramirez.2012, lin.2020b}.

\Ac{OC} requires accurate models of the \ac{GH} building environment and the crop growth. Crop models can offer insights into how various factors influence crop growth and yield. Existing \ac{OC} approaches \cite{vanStraten.2013, engler.2021} use models customized for specific crops and problems, which necessitates extensive experimentation for each adaptation and limits their general applicability. 
In contrast, numerous models have been developed for the general description of plants in agriculture and biology. Among them are several tomato models developed in the last 30 years \cite{jones.1991, heuvelink.1999, vanStraten.2010, SEGINER.1994, marcelis.2008, MartnezRuz.2019, Boote.2012}, indicating the continuous interest of research and industry in more accurate models for various applications.
However, control engineering requires a mathematically analytical description of models and most of these models were not designed for that use-case and are not given as a set of equations but as a defined combination of look-up tables. 
Hence, they are often implemented as a pure simulation black box model with no interfaces to the inside of the models, and the models themselves frequently include discontinuities, and only guarantee good performance within unknown boundaries. As a result, selecting a model that is accurate, computationally efficient, and appropriate for the use in \ac{OC} remains a challenge.

Previous studies have used different approaches to select or synthesize crop models for control engineering, including a model comparison \cite{ramirez.2003} of the \textit{TomSim} \cite{heuvelink.1999} and reduced \textit{TOMGRO} model \cite{jones.1999} and combining components from existing models \cite{tap.2000, vanOoteghem.2007, kuijpers.2019}. 
While previous studies have used a common model structure to combine model components into new models \cite{kuijpers.2019}, we do not use this approach since the individual components lack modularity and explainability. Instead we compare existing models regarding their performance and applicability in \ac{OC}. To the best of our knowledge, there is a lack of comparative studies on tomato crop models that include performance metrics regarding \ac{OC} results, computational time, and applicability.

In this paper, we compare three established crop models. The \textit{SIMPLE} crop model by Zhao et al. \cite{Zhao.2019} is a simple but generic crop model that can be used for various crops. The reduced \textit{TOMGRO} model by Jones et al. \cite{jones.1999} is designed for tomatoes only, and agricultural scientists and farmers commonly use the tomato model of the \textit{DSSAT} toolbox \cite{Jones.2003}. The comparison evaluated the model's performance when applied to \ac{OC} approaches in a \ac{GH} environment. We addressed the remaining challenges that may impact the implementation of these models in real-world agricultural settings.

Our key contributions can be summarized as follows:
\begin{itemize}
\item Comparison of the accuracy of \textit{SIMPLE}, reduced \textit{TOMGRO}, and \textit{DSSAT},  using a \ac{GH} data set from the Autonomous \ac{GH} Challenge in the Netherlands \cite{hemming.2020a,Hemming.2020b}.
\item Comparison of the structure and applicability of these models in \ac{OC} and re-formulation of two models as state-space models to be used in control theory. Combination of the models with a \ac{GH} environment model.
\item Incorporation of the \ac{GH}-crop-model into an \ac{OC} approach, and analysis of the results for their performance and meaning with recommendations for future users.
\end{itemize}

The paper is organized as follows. Sec. \ref{sec:models} provides an overview of the tomato models and the \ac{GH} environment. In Sec. \ref{sec:optimalcontrol-framwork}, we introduce the \ac{OC} approach. Sec. \ref{sec:comparison} discusses the validation of the models and presents the comparison results. The ensuing Sec. \ref{sec:discussion} delves into the implications of our findings and explores the limitations of this study and future research directions. Finally, in Sec. \ref{sec:conclusion}, we conclude the paper.

\section{Models for the Growth of Tomatoes}
\label{sec:models}
Before a comparison can be conducted, an introduction to the three crop models and the \ac{GH} model is given.
\subsection{SIMPLE Model}
The parameters of the \textit{SIMPLE} model \cite{Zhao.2019} have been calibrated carefully using a large arable farming experimental data set for various crops. The model  represents up to 14 crops, including tomatoes, and 22 cultivars with the modification of just 13 crop parameters, of which four are cultivar-specific and nine are species-specific. Although the model takes into account various factors like plant phenology, the impact of photosynthesis on growth, the influence of CO$_{2}$ concentration, drought stress, and radiation interception, it has some limitations, including the exclusion of vernalization effects and the lack of nutrient dynamics. 
The state space-model for  day $i$ can be described as
\begin{equation} 
    \bs{x}_{\text{s}, i+1} =  \bs{f}_\text{s}(\bs{x}_{\text{s}, i}, \bs{u}_{\text{s}, i})
    \label{eqn:statesimple}
\end{equation}
with the state of the  \textit{SIMPLE} model
\begin{equation} 
    \bs{x}_{\text{s}, i}= \begin{bmatrix} m_{B,i}  &  \tau_{i}  &  I_{50B,i} \end{bmatrix}^\mathsf{T} 
\end{equation}
where $m_{B,i}$ is the tomato biomass in \si{\kilogram\per\metre\squared}, $\tau_{i}$ is the cumulative temperature in \si{\degreeCelsius\day} (temperature integrated over days) and $I_{50B,i}$ is a value for the leaf senescence on day $i$ in \si{\degreeCelsius\day} \cite{Daniels.2023}.	
The yield for the tomatoes is obtained with $m_{\text{fruit}} = H\!I \ m_{B, N}$, where $H\!I$ is the harvest index ($H\!I=0.68$ for tomatoes) and $m_{B, N}$ is the biomass of the plant on the last day.

Following the method employed in \cite{Daniels.2023}, we assume that the temperature remains constant throughout the day. However, we additionally treat the CO\textsubscript{2} concentration as a controllable variable. The input vector for this model for day $i$ is
\begin{equation} \label{eq:usimple}
    \bs{u}_{\text{s}, i}= \begin{bmatrix}T_{i}  &  D_{i}  &  R_{i} & C_{CO2,i}  \end{bmatrix}^\mathsf{T} 
\end{equation}
where $T_{i}$ is the mean temperature in \textcelsius, $D_{i}$ is the the relative level of drought between $[0,1]$, $R_{i}$ is the solar radiation in \si{\mega\joule\per\metre\squared\per\day} and $C_{CO2,i}$ is the CO\textsubscript{2} concentration in \si{ppm}.
\subsection{Reduced TOMGRO Model}
The reduced  \textit{TOMGRO} model \cite{jones.1999} is a simplified version of the original, widely used \textit{TOMGRO} model \cite{jones.1991}. Compared to the original model, where the state dimension is 69, the reduced \textit{TOMGRO} model contains a state with only a dimension of 5 and an input dimension of 3. This model was developed under \ac{GH} conditions.

Similar to \eqref{eqn:statesimple}, we propose the reduced \textit{TOMGRO} model describing the state of the next day $i+1$  as
 \begin{equation} 
    \bs{x}_{\text{t}, i+1} = \bs{f}_\text{t}(\bs{x}_{\text{t}, i}, \bs{u}_{\text{t}, i}) \, .
    \label{eqn:statetomgro}
\end{equation}

The state vector of this model for day $i$ is given as
\begin{equation} \label{eqn:stateTOMGRO}
    \bs{x}_{\text{t}, i}= \begin{bmatrix}N_{i}  &  L\!A\!I_{i}  &  W_{i} & W_{f,i} & W_{m,i}  \end{bmatrix}^\mathsf{T}  ,
\end{equation} 
where $N_{i}$ is the number of main stem nodes, $L\!A\!I_{i}$ is the leaf area index in  \si{\metre\squared\per\metre\squared} (ratio of leaf area per ground area), $W_{i}$ is the total plant weight in  \si{\kilogram\per\metre\squared}, $ W_{f,i}$ is the fruit dry weight in \si{\kilogram\per\metre\squared}, and $W_{m,i}$ is the mature fruit dry weight in  \si{\kilogram\per\metre\squared}.
 
The input vector of this model for day $i$ is given as
\begin{equation} \label{eq:utomgro}
    \bs{u}_{\text{t}, i}= \begin{bmatrix}T_{i} & T_{d,i} & R_{i} & C_{CO_{2},i} \,  \end{bmatrix}^\mathsf{T} ,
\end{equation}
where $T_{i}$ is the average temperature of the whole day in \textcelsius, $T_{d,i}$ is the average temperature of daytime in \textcelsius, $R_{i}$ is the solar radiation on day $i$ in \si{\mega\joule\per\metre\squared\per\day}, 
$C_{CO_{2},i}$ is the concentration of CO\textsubscript{2} in ppm. 
From $R_{i}$, we obtain the \acl{PPFD} $P\!P\!F\!D_{i}= R_{i}/\SI{0.037}{\mega\joule\second\per\micro\mole\per\day}$ for white light. 

As the original model description contains piece-wise defined functions and is therefore not always differentiable, we use the smoothing function for maximum operators (cf. \cite{Daniels.2023}), as well as a smoothing function for Heaviside step function $H(x) \approx H_\epsilon (x) =  \displaystyle{\frac{1}{1+e^{-\epsilon x}}}$
to allow for gradient-based optimization approaches.
 For a big positive $\epsilon$, the smoothing function converges to the standard Heaviside step function. 
 
\subsection{DSSAT}
The \textit{DSSAT} toolbox \cite{Jones.2003} contains widely used and well-established models in arable agriculture and has been used for decision support processes for a long time. For tomatoes it uses the \textit{CROPGRO} model \cite{Boote.2012}. However, its implementation is only available as a complete Windows application. It requires manual file uploads, making it unsuitable for integration into \ac{OC} processes. Therefore, we will use this model only as a well-established comparison.
\subsection{Greenhouse Model}

The \ac{GH} model described in \cite{vanStraten.2010} is employed for both the \textit{SIMPLE} crop model \cite{Zhao.2019} and the reduced \textit{TOMGRO} model \cite{jones.1999}. The \ac{GH} model comprises the primary \ac{GH} compartment, which is linked to three subsystems. These subsystems serve as the three control inputs for the \ac{GH}, which are the heat supply through a heating pipe system, the supply of CO\textsubscript{2}, and the window ventilation system for air flow. For a more comprehensive explanation of the \ac{GH} structure, see \cite{vanStraten.2010}.

The state-space equation for the GH model describes the state of the \ac{GH} on an hourly basis and is given as
\begin{equation} 
 \label{eqn:state}
             \bs{x}_{\text{gh}, i+1}=  \bs{f}_\text{gh}( \bs{x}_{\text{gh}, i},  \bs{u}_{\text{gh}, i}) \, .
    \end{equation}
The state vector  of the \ac{GH} model is
\begin{equation} 
    \bs{x}_{\text{gh}, i}= \begin{bmatrix}T_{g,i}  &  T_{s,i}  & T_{p,i} & C_{CO_{2},i}  & C_{H_{2}O,i} \end{bmatrix}^\mathsf{T} \, ,         
\end{equation}
which are the greenhouse temperature $T_{g,i}$ in \si{\degreeCelsius}, soil temperature $T_{s,i}$ in \si{\degreeCelsius}, pipe temperature $T_{p,i}$ in \si{\degreeCelsius}, CO\textsubscript{2} concentration $C_{CO_{2},i}$ in \si{ppm}, and the water vapor concentration $C_{H_{2}O,i}$ in \si{\kilogram\per\cubic\metre} for the $i^\text{th}$ hour.
The three control inputs of the \ac{GH}  
    \begin{equation} \label{eq:inputgreenhouse}
            \bs{u}_{\text{gh},i}= \begin{bmatrix}u^{vp}_{q,i}  & u^{Ap}_{v,i} & u^{vp}_{CO_2,i}  \end{bmatrix}^\mathsf{T} 
        \end{equation}
are directly connected to the main \ac{GH} compartment. The underlying systems of the control inputs are simplified. The heat and CO\textsubscript{2} supply are expressed as a valve position $u^{vp}_{q,i}$ and $u^{vp}_{CO_2,i}$  which quantifies the supply, respectively. The ventilation system is expressed by the aperture of the windows from the windward and lee side of the greenhouse compartment. 
Both windows are considered for simplicity through one control input $u^{Ap}_{v,i}$. 

Note that the the valve positions $u^{vp}_{q,i}$  and $u^{vp}_{CO_2,i}$ are between $[0,1]$, where 0 indicates that there is no supply and the valve is closed. 
The control input $u^{Ap}_{v,i}$ lies in the interval $[0,2]$ since we consider both the windward and lee side of the \ac{GH} compartment.  
Furthermore, the \ac{GH} compartment is subject to environmental influences as it is a partially open system. 
Thus, the overall dynamics of the \ac{GH} model also include six parameters, which are summarized in the external disturbance vector
	 	\begin{equation} 
	 	\label{eqn:externalinput}
            d_{g,i}= \begin{bmatrix}
            R_{\text{out},i} \ \, T_{\text{out},i} \ \, v_i \ \, T_{\text{soil},i} \ \, C_{{H_{2}O\,\text{out},i}} \ \, C_{CO_{2}\,\text{out},i}  \end{bmatrix}^\mathsf{T}.
        \end{equation} 
with the solar radiation $R_{\text{out},i}$ in \si{\watt\per\square\metre}, the temperature outside the greenhouse $T_{\text{out},i}$ in \si{\degreeCelsius}, wind speed $v_i$ in \si{\metre\per\second}, temperature of the subsoil $T_{\text{soil},i}$ in \si{\degreeCelsius}, the humidity $C_{{H_{2}O\,\text{out},i}}$ in \si{\kilogram\per\cubic\metre} and concentration of CO\textsubscript{2}  $C_{CO_{2}\,\text{out},i}$ in \si{ppm} outside the \ac{GH}. 
Just as with the crop models, smoothing functions are necessary to handle discontinuities. 
For example, an approximation of the absolute value is ${|x| \approx \sqrt{x^2 + \mu}}$ where $\mu$ is a small positive constant \cite{ramirez2014}. The derivative of the approximation exists for all  $x\in\mathbb{R}$. 

\subsection{Integration of the models}
\label{sec:integration}
In a comparison of the \ac{OC} results of the two harvesting models, the respective state-space model is used in combination with the \ac{GH} model. 
Since the \ac{GH} model operates on an hourly basis while the crop models are evaluated once a day, evaluating the combined model is not a simple task. Consequently, the crop model receives the mean value of the \ac{GH} conditions as input, and the combined state vector comprises the state of the tomato model along with the hourly state of the greenhouse, i.e., 
\begin{align}\label{eq:combineStates}
    \bs{x}_i = \begin{bmatrix} \bs{x}_{\text{s/t},i}^\top & \bs{x}_{\text{gh},24i}^\top &\bs{x}_{\text{gh},24i+1}^\top &.. & \bs{x}_{\text{gh},24i+23}^\top\end{bmatrix}^\top,
\end{align}
where $\bs{x}_{\text{s/t},i}$ is either the  state resulting from the \textit{SIMPLE} model $\bs{x}_{\text{s},i}$ or the state based on the \textit{TOMGRO} model $\bs{x}_{\text{t},i}$. 

In order to account for and optimize the performance of the entire system comprising the tomato crop models within a greenhouse setting, we amalgamate them into a single model
\begin{align}
    \bs{x}_{i+1}= \bs{f}(\bs{x}_{i}, \bs{u}_{\text{gh}, i}) \, ,
    \label{eq:combined}
\end{align}
where the components of $\bs{x}_{i+1}$ are computed by
\begin{align}
    \bs{x}_{\text{s/t}, i+1} &=  \bs{f}_\text{s/t}\left(\bs{x}_{\text{s/t}, i},\bs{g}_\text{s/t}\left(  {\textstyle\frac{1}{24}\sum_{j=0}^{23} \bs{x}_{\text{gh}, i+j}} \right) \right) \\
    \bs{x}_{\text{gh}, i+j+1} &= \bs{f}_\text{gh}( \bs{x}_{\text{gh}, i+j},  \bs{u}_{\text{gh}, i}) \quad \forall j \in [0,23].
\end{align}
and arranged in the same form as in \eqref{eq:combineStates}. The function $\bs{g}_\text{s/t}\left( \cdot \right) $ maps the mean states of the \ac{GH} to the inputs for either the \textit{SIMPLE} or the \textit{TOMGRO} model.
In the following section, the combined model \eqref{eq:combined} of  tomato growth and \ac{GH} is used in an \ac{OCP}.

\section{Optimal Control for Plant Growth}
\label{sec:optimalcontrol-framwork}
In the following, the framework for the \ac{OC} of plant growth is presented. 
The goal is to find a sequence of  input values for the greenhouse \eqref{eq:inputgreenhouse} resulting in a maximum economic yield. The inputs of the greenhouse over the optimization period $N$ is combined to the input sequence vector ${\bs{U} = \left[ \bs{u}_{\text{gh},0}^\top, \bs{u}_{\text{gh},1}^\top, ..., \bs{u}_{\text{gh},N-1}^\top\right]^\top}$.

The optimization selects one from all possible input sequences $\bs{U}$ that minimizes the cost function
\begin{align}\label{eq:cost}
    J(\bs{U},\bs{x}_N) = \sum_{i=0}^{N-1}  l( \bs{u}_{\text{gh},i}) - V(\bs{x}_N),
\end{align}
while the input and state sequences are a solution to the system model \eqref{eq:combined} that consists of the greenhouse \eqref{eqn:state} and the tomato model, either the \textit{SIMPLE} model \eqref{eqn:statesimple} or the \textit{TOMGRO} model~\eqref{eqn:statetomgro}. 
It is evaluated over a given growth period of $N$ days. 
The cost function \eqref{eq:cost} represents the negative economic yield, thus the terminal cost is used with a negative sign, similar to \cite{Daniels.2023}. 
The stage cost $l( \bs{u}_{\text{gh},i})$ represents the energy cost used for the growing process and is given as 
$l( \bs{u}_{\text{gh},i}) = \bs{r}^\mathsf{T} \bs{u}_{\text{gh},i}$,
where 
$\bs{r} \in \mathbb{R}^{3}$ is a  weight vector. The terminal cost $V(\bs{x}_N)$ gives the yield at the harvest on day $N$. The price evolves linearly, leading to a linear cost term, i.e.  $V(\bs{x}_N) =\bs{q}^\mathsf{T} \bs{x}_N$, with the weight $\bs{q}\in\mathbb{R}^{3}$. 
The temperature, ventilation and CO\textsubscript{2}  can be controlled in a \ac{GH} within lower and upper bounds. Therefore, we add constraints $\bs{u}_{\text{gh},i}\in \mathcal{U}$ to the \ac{OCP}, yielding 
\begin{subequations}\label{eq:opt}
    \begin{align}
        \bs{U}^* =& \argmin_{\bs{U}}  J(\bs{U},\bs{x}_N)  \\
        \text{s.t. }& \bs{x}_{i+1}= \bs{f}(\bs{x}_{k},\bs{u}_{\text{gh},i})  &\forall i \in[0,N-1] \\
        & \bs{u}_{i}\in \mathcal{U}  \quad &\forall i \in[0,N-1] \\
        & \bs{x}_0 = \bs{x}_\text{init},
    \end{align}
\end{subequations}
where $\bs{U}^*$ is the sequence of optimal inputs of the \ac{GH}. The states on day $i = 0$ are given as initial states $\bs{x}_\text{init}$.  
The optimal input sequence $\bs{U}^*$ results in a state sequence that yields the smallest value of the cost function \eqref{eq:cost}. 
The OCP is implemented in Python and solved with CasADi \cite{Andersson.2019}. 

\section{Comparison of the crop models}
After the defintion of the models, we now compare these models regarding their structure and their validity compared to the \textit{DSSAT} model and a \ac{GH} data set. Additionally, we compare the findings of the \ac{OC} study.
\label{sec:comparison}
\subsection{Structure}
\label{sec:structure}
As already described in Sec. \ref{sec:integration}, the two tomato models follow a similar state-space representation, which are integrated with the \ac{GH}.
Furthermore, their inputs,
denoted as $\bs{u}_{\text{s}, i}$ and $\bs{u}_{\text{t}, i} $, 
have a similar structure and contain variables such as daily average temperature, solar radiation, and CO\textsubscript{2}. However, there are some differences between the two models in terms of inputs. For instance, the \textit{SIMPLE}  model includes an input for a drought factor $D$, which is set to 0 in a \ac{GH} environment under the assumption that we guarantee sufficient water supply for each day. 
Another key difference is that \textit{TOMGRO} follows a more complex model approach, e.g., it has a state with five elements \eqref{eqn:stateTOMGRO}, whereas the state of the \textit{SIMPLE}  model only considers three elements \eqref{eqn:statesimple}. While the \textit{TOMGRO} model state is readily accessed and measurable through simple techniques such as counting stem nodes, assessing leaf area, and weighing different parts of the plant, the state of the \textit{SIMPLE} model is not as easily observable. It requires expert knowledge along with analytical methods. Nonetheless, both models share a common state, the dry weight of the fruit ($m_{\text{fruit}}$ for the \textit{SIMPLE} model and $W_f$ for the \textit{TOMGRO} model), which is of interest in this study as it helps in obtaining the optimal yield.
The results are given in fresh weight to compare it with experimental data.

\subsection{Validation of Models}
\label{sec:validation}
For the validation of both models, we use as real-world ground truth a data set consisting of \ac{GH} climate data compiled from different experiments \cite{Hemming.2020b}.
This data set, which we compare to both models, is generated from the experiments of five different multi-disciplinary teams.
As an example and for brevity, we only show the comparison for one set (The Automators team) for our validation. We use the environmental conditions of the experiment, including temperature, solar radiation, and CO\textsubscript{2} concentration, for each day.

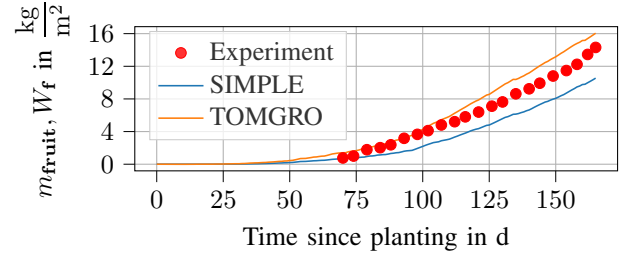
\begin{figure}[ht]
	\centering 
	\vspace{-0.5cm}
\begin{tikzpicture}

\definecolor{darkgray176}{RGB}{176,176,176}
\definecolor{darkorange25512714}{RGB}{255,127,14}
\definecolor{lightgray204}{RGB}{204,204,204}
\definecolor{steelblue31119180}{RGB}{31,119,180}

\begin{axis}[
height=3.5cm,
legend cell align={left},
legend style={
  fill opacity=0.8,
  draw opacity=1,
  text opacity=1,
  at={(0.03,0.97)},
  anchor=north west,
  draw=lightgray204
},
minor xtick={},
minor ytick={},
tick align=outside,
tick pos=left,
width=8cm,
x grid style={darkgray176},
xlabel={Time since planting in d},
xmajorgrids,
xmin=-8.25, xmax=173.25,
xtick style={color=black},
xtick={-25,0,25,50,75,100,125,150,175},
y grid style={darkgray176},
ylabel={\(\displaystyle m_\mathbf{fruit},W_\mathbf{f}\) in \(\displaystyle \mathrm{\frac{kg}{m^2}}\)},
ymajorgrids,
ymin=-0.812918417301494, ymax=16.8561853558258,
ytick style={color=black},
ytick={0,4,8,12,16}
]
\addplot [draw=red, fill=red, mark=*, only marks]
table{%
x  y
70 0.767
74 0.999
79 1.777
84 2.025
88 2.379
93 3.174
98 3.66
102 4.114
107 4.811
112 5.209
116 5.804
121 6.392
126 7.112
130 7.631
135 8.626
140 9.239
144 9.927
149 10.807
154 11.477
158 12.239
162 13.452
165 14.311
};
\addlegendentry{Experiment}
\addplot [semithick, steelblue31119180]
table {%
0 0
1 2.43674068413549e-05
2 6.03502808098544e-05
3 0.000151267458389586
4 0.000203911125475844
5 0.000233347031868218
6 0.000282560241067372
7 0.000312069228921196
8 0.000409609487855019
9 0.000477403104615676
10 0.000619640895671799
11 0.000694153474196344
12 0.00108402857481072
13 0.0014218469547255
14 0.00170112864331437
15 0.00235087795521694
16 0.00298791004710111
17 0.0032095119051274
18 0.00337953717326475
19 0.00357535963262975
20 0.00416797025123469
21 0.00444409261078599
22 0.005951024672374
23 0.00709676862236841
24 0.00781045432793006
25 0.00821210776451261
26 0.0101886572454236
27 0.0124737784090221
28 0.0132543374383
29 0.016997085506319
30 0.0181774564486438
31 0.0191988399361757
32 0.0265389814982133
33 0.028896976571083
34 0.0363450469673398
35 0.0461479488575086
36 0.0517700191823603
37 0.0572342086973954
38 0.0606248652898861
39 0.0709674470466897
40 0.0797928569160651
41 0.0847118724392604
42 0.100274627271173
43 0.109693608414207
44 0.120436003016953
45 0.134094054314973
46 0.149219752476043
47 0.155864080846883
48 0.174846414838834
49 0.180711162223755
50 0.205018529594832
51 0.216558028500006
52 0.250777091533802
53 0.28699934600997
54 0.329531974700792
55 0.337544428915428
56 0.344654147755318
57 0.369520711368574
58 0.391643648947253
59 0.419500242292492
60 0.429831626616507
61 0.454609472653538
62 0.474856440981298
63 0.484661896839278
64 0.517629647992813
65 0.557238397477335
66 0.597295162623881
67 0.619099879058623
68 0.641729864824995
69 0.657048604651462
70 0.668828210770767
71 0.676341729845902
72 0.727891641239033
73 0.753034421541875
74 0.765172569412609
75 0.791569012233188
76 0.822593825083119
77 0.861352358556318
78 0.874628087513993
79 0.945112486146686
80 1.00176789299803
81 1.01277867171383
82 1.0499838140406
83 1.11766468120817
84 1.13744760492337
85 1.17481280331572
86 1.18432179697341
87 1.24061140217189
88 1.31756938555826
89 1.34832903330872
90 1.40307593931445
91 1.42616238722672
92 1.50353150573198
93 1.59162057942115
94 1.64132331867488
95 1.66497344500885
96 1.69647224992486
97 1.81546986351622
98 1.93712308233125
99 2.06030233636178
100 2.18515816446575
101 2.31344854925557
102 2.42892957339003
103 2.53435921917687
104 2.64623775423392
105 2.69798629882699
106 2.74973496819342
107 2.87548104230005
108 2.93752285761231
109 2.99861268263072
110 3.07065892363313
111 3.19455179691099
112 3.32685425319012
113 3.44054971432828
114 3.5738414238115
115 3.69396071898128
116 3.81036517563375
117 3.94204191791059
118 4.08218022817321
119 4.19134140774809
120 4.28950826971027
121 4.35969407655679
122 4.49232377722395
123 4.62636060168339
124 4.75784552979738
125 4.80289616889564
126 4.95310367866475
127 5.10735194936238
128 5.26486601548505
129 5.4256939450653
130 5.58392030135227
131 5.73985857598314
132 5.8254148812029
133 5.96541005977589
134 6.11659654519058
135 6.136826758662
136 6.20791139349992
137 6.32121661143441
138 6.41186728168412
139 6.5517221835682
140 6.65060179416265
141 6.81167701570655
142 6.98891951973826
143 7.15548658967606
144 7.31366079450649
145 7.46169004550142
146 7.61553908462608
147 7.70154606049421
148 7.85253117654042
149 7.97246668099681
150 8.07395383514306
151 8.22625484975063
152 8.39158033917889
153 8.52694473606948
154 8.71512866888351
155 8.90447595414035
156 9.09406816207552
157 9.26563404751495
158 9.46989429136122
159 9.58220080025232
160 9.76074388777123
161 9.79887069108814
162 9.98307700216136
163 10.163111759337
164 10.3488165948409
165 10.5475352583197
};
\addlegendentry{SIMPLE}
\addplot [semithick, darkorange25512714]
table {%
0 0.001904
1 0.000122790728424587
2 -0.00114608422339483
3 -0.00198802578858541
4 -0.00296435886291252
5 -0.00444585051432977
6 -0.00575468991590465
7 -0.00727723663841007
8 -0.00786045739407817
9 -0.00884054535293039
10 -0.00887646461874998
11 -0.00977733670479652
12 -0.00709552331740964
13 -0.00518838624013203
14 -0.00382580747471724
15 0.00108079410738362
16 0.00525408683794693
17 0.00560818605553364
18 0.0052405659192954
19 0.00488019965255715
20 0.00727065732603048
21 0.00703244489722756
22 0.015159136026752
23 0.0204466163107367
24 0.0224336664299775
25 0.0222926915753446
26 0.0305660198834682
27 0.0397807422608533
28 0.0410105040523011
29 0.0547690471878057
30 0.0568332286450695
31 0.0580466864270657
32 0.0811918892936316
33 0.0864227183707483
34 0.107867731902484
35 0.134861094238373
36 0.148042848698634
37 0.160692276186621
38 0.166635694377575
39 0.189069935900312
40 0.207214268699845
41 0.215149651745779
42 0.248067591802539
43 0.26573991591869
44 0.287039000545023
45 0.314346423559454
46 0.343868917699687
47 0.353558402442646
48 0.391166368555167
49 0.398885579997245
50 0.448315710816147
51 0.469921662145674
52 0.538035707753507
53 0.611017976407977
54 0.692653926204596
55 0.702807487551495
56 0.711406100389686
57 0.76487592540584
58 0.812687960359275
59 0.8742909252663
60 0.890527836534649
61 0.943163317622027
62 0.98422062190392
63 0.998541604443001
64 1.07467633746673
65 1.16566983268095
66 1.25795238050488
67 1.30525284018567
68 1.35511345714695
69 1.3856633621882
70 1.40495403579464
71 1.41120627963795
72 1.51906293015452
73 1.57971557949033
74 1.60081519977774
75 1.66397952294204
76 1.74052327576276
77 1.83918323021925
78 1.86099105733552
79 1.9899535033631
80 2.10719072591283
81 2.11995274323459
82 2.21892941789512
83 2.34818321995972
84 2.38914279369497
85 2.48524434645651
86 2.49031226926371
87 2.61056274920318
88 2.75250036548702
89 2.82912297240717
90 2.94912626473338
91 2.99547774064896
92 3.13854608946114
93 3.29221486174185
94 3.40713775915032
95 3.45259149651104
96 3.52903058679155
97 3.71338975754233
98 3.89992854575029
99 4.08734510370251
100 4.27527045524561
101 4.46627319368884
102 4.64991229614715
103 4.82584276414096
104 5.00892526579459
105 5.13180259981864
106 5.25468851710149
107 5.45213949577475
108 5.58164663153868
109 5.70984047715372
110 5.85189174803622
111 6.04689417780348
112 6.24685808044825
113 6.43066494935044
114 6.62895592568334
115 6.81402571885804
116 6.99602136011907
117 7.19411045351977
118 7.39587892146967
119 7.56956565761541
120 7.73971843374209
121 7.87527193938773
122 8.07661232857511
123 8.27356046993519
124 8.47148415064938
125 8.55683706211571
126 8.76972317145659
127 8.98535121382689
128 9.20107455419322
129 9.41733535873167
130 9.62928667348119
131 9.84262613917674
132 9.98793451035247
133 10.1875848109744
134 10.3944801846482
135 10.3861536179146
136 10.5109557348139
137 10.6831417069337
138 10.8309761905013
139 11.0278120003017
140 11.1831799989587
141 11.3971138218625
142 11.627049012328
143 11.8457285135195
144 12.0564371603205
145 12.2525600094617
146 12.4572606995244
147 12.6021821550776
148 12.8138366773771
149 12.9939346149407
150 13.1610814600876
151 13.3750351793195
152 13.5919375756726
153 13.7778340040051
154 13.9978646967258
155 14.2153511927708
156 14.4309970878597
157 14.629927115932
158 14.8004633762487
159 14.9468575088041
160 15.1566790554305
161 15.1797959053844
162 15.3963801815152
163 15.6074140558094
164 15.8282209271188
165 16.0530442752291
};
\addlegendentry{TOMGRO}
\end{axis}

\end{tikzpicture}
	\caption{Comparison of the tomato fresh weight computed with the \textit{SIMPLE}  and \textit{TOMGRO} model with experimental data.} 
	\vspace{-0.2cm}
	\label{fig:validation}
\end{figure}

In Fig. \ref{fig:validation}, we observe that both models present a similar behavior for around $160$ days of cultivation in comparison to the experimental data. 
However, the \textit{TOMGRO} model fits the data better and yields about \SI{2}{\kilogram\per\metre\squared} more biomass than the experimental data, while the \textit{SIMPLE} model results in an underestimate of \SI{4}{\kilogram\per\metre\squared}, under the same conditions.

\begin{figure}[ht]
	\centering 
	\vspace{-0.5cm}
\begin{tikzpicture}

\definecolor{darkgray176}{RGB}{176,176,176}
\definecolor{darkorange25512714}{RGB}{255,127,14}
\definecolor{lightgray204}{RGB}{204,204,204}
\definecolor{steelblue31119180}{RGB}{31,119,180}

\begin{axis}[
height=3.5cm,
legend cell align={left},
legend style={
  fill opacity=0.8,
  draw opacity=1,
  text opacity=1,
  at={(0.03,0.97)},
  anchor=north west,
  draw=lightgray204
},
minor xtick={},
minor ytick={},
tick align=outside,
tick pos=left,
width=8cm,
x grid style={darkgray176},
xlabel={Time since planting [d]},
xmajorgrids,
xmin=-6.05, xmax=127.05,
xtick style={color=black},
xtick={-20,0,20,40,60,80,100,120,140},
y grid style={darkgray176},
ylabel={\(\displaystyle m_\mathbf{fruit},W_\mathbf{f}\) in \(\displaystyle {\mathrm{\frac{kg}{m^2}}}\)},
ymajorgrids,
ymin=-0.670015, ymax=14.070315,
ytick style={color=black},
ytick={0,3,6,9,12}
]
\addplot [semithick, red]
table {%
28 0
29 0
30 0.0003
31 0.0012
32 0.0029
33 0.0057
34 0.0097
35 0.0154
36 0.0226
37 0.0322
38 0.0431
39 0.0638
40 0.0892
41 0.1221
42 0.1631
43 0.2181
44 0.2884
45 0.3816
46 0.4956
47 0.6132
48 0.7584
49 0.9407
50 1.1209
51 1.3169
52 1.5275
53 1.7272
54 1.916
55 2.0923
56 2.2675
57 2.4539
58 2.6593
59 2.8327
60 2.99
61 3.2088
62 3.466
63 3.7352
64 3.9882
65 4.1747
66 4.3522
67 4.62
68 4.892
69 5.1657
70 5.3682
71 5.598
72 5.8537
73 6.1521
74 6.3949
75 6.6735
76 6.9764
77 7.2756
78 7.5691
79 7.8564
80 8.0957
81 8.2898
82 8.5307
83 8.8169
84 9.0804
85 9.3297
86 9.5827
87 9.8051
88 9.9938
89 10.1938
90 10.3651
91 10.5504
92 10.7402
93 10.881
94 11.0542
95 11.2294
96 11.3815
97 11.5567
98 11.7011
99 11.8414
100 11.9909
101 12.1509
102 12.2793
103 12.3836
104 12.4872
105 12.5523
106 12.6454
107 12.7905
108 12.8669
109 12.9533
110 13.0427
111 13.1128
112 13.1729
113 13.2225
114 13.2514
115 13.2678
116 13.2934
117 13.3003
118 13.3693
119 13.3937
120 13.4003
121 13.4003
};
\addlegendentry{DSSAT}
\addplot [semithick, steelblue31119180]
table {%
0 0
1 0.000297350971489009
2 0.000826224293510826
3 0.00148592413948614
4 0.00172942400940233
5 0.00196787819848889
6 0.00249720740400547
7 0.00313282498977328
8 0.00420872706322755
9 0.00513837594434572
10 0.00606574476510766
11 0.00729694869767827
12 0.00860092953624112
13 0.0103509896433377
14 0.0128673556146345
15 0.0149810384493339
16 0.0169970567443408
17 0.0176667760500045
18 0.0235932865648671
19 0.0311707354020019
20 0.0358822540844783
21 0.0395591952346023
22 0.0452850647835127
23 0.0526171442246644
24 0.065009897071956
25 0.0749403678829347
26 0.0771757195855289
27 0.0868180468582321
28 0.0946738283119357
29 0.106055849904909
30 0.111829018335392
31 0.133822472922829
32 0.140191708528269
33 0.144047850270025
34 0.170989531142527
35 0.178253609429731
36 0.216199216212405
37 0.25755430199409
38 0.310756851326553
39 0.373332569155059
40 0.433434559746346
41 0.49164394227571
42 0.576724036342559
43 0.672815235907289
44 0.784589496186441
45 0.892990574786438
46 0.943348585700833
47 0.964376308105957
48 0.985289044760634
49 1.09721218642686
50 1.20892578648977
51 1.33150252876721
52 1.43412259304666
53 1.56772139732597
54 1.72553160765449
55 1.87277649304073
56 2.02704028383816
57 2.19682313933063
58 2.36937962947762
59 2.53992571924997
60 2.73165974410981
61 2.88712036871329
62 3.04172296445352
63 3.24059179661753
64 3.43911677102924
65 3.63300604440701
66 3.82030482172834
67 3.97364848399876
68 4.14007767506399
69 4.28094293537867
70 4.46186713901258
71 4.62444063085366
72 4.84469503911248
73 4.89866456199688
74 5.09565573612584
75 5.30331734055044
76 5.54871864828047
77 5.74843585212785
78 5.94633434337861
79 6.16291012530814
80 6.3772854332769
81 6.62644261431483
82 6.847529848684
83 7.04127863078673
84 7.24746794926904
85 7.4710283432899
86 7.69490037622454
87 7.92296819073123
88 8.02283571812265
89 8.10802531293073
90 8.19445885325977
91 8.27944687773481
92 8.36337664972476
93 8.44629435451413
94 8.52934872397397
95 8.64572280696661
96 8.84734274483236
97 9.00975944898566
98 9.09976619064581
99 9.16137685359692
100 9.22193825325704
101 9.38715654969535
102 9.54218574815017
103 9.71221272700829
104 9.81404203455219
105 10.0103167680565
106 10.2458460772829
107 10.4838714125215
108 10.6763480434805
109 10.8095447316302
110 11.0275709601872
111 11.2036562787
112 11.4250161451868
113 11.6415653928141
114 11.791644143306
115 11.8644889026752
116 12.0932270568255
117 12.3184985829593
118 12.5053585762744
119 12.731371829542
};
\addlegendentry{SIMPLE}
\addplot [semithick, darkorange25512714]
table {%
0 0.0028
1 0.00193024256398179
2 0.0011353733614063
3 0.00106624386164328
4 0.00384343876787485
5 0.00729505843943589
6 0.0103117765339627
7 0.0125159725506181
8 0.016654494534319
9 0.0214381056740046
10 0.0256830175845029
11 0.0335063264500062
12 0.0437594508502617
13 0.0542894036834638
14 0.0695514479807671
15 0.0893445961088597
16 0.111847671141043
17 0.118754810731741
18 0.142518066655991
19 0.169681641540984
20 0.192531190121398
21 0.218937217196959
22 0.257436334605912
23 0.297667738942155
24 0.339738639836315
25 0.376221285564042
26 0.386249117852944
27 0.431448465994798
28 0.459439719797074
29 0.508499420675533
30 0.533390962309323
31 0.591879550512407
32 0.619258351410987
33 0.638120732756471
34 0.699706405022204
35 0.719890515993619
36 0.801975833219086
37 0.882151308112636
38 0.976105121953847
39 1.07298175582685
40 1.16954361218992
41 1.26089922118345
42 1.36614126491784
43 1.47169309752541
44 1.57714937619572
45 1.68136132938436
46 1.74003356249101
47 1.77088215742808
48 1.81320933285798
49 1.90440363055595
50 1.99698946243542
51 2.09898411839941
52 2.20911024753608
53 2.33276500624766
54 2.46271046004467
55 2.59489346982195
56 2.7305985874745
57 2.87264888626146
58 3.01732297696808
59 3.16544391197604
60 3.32128147544892
61 3.47526465525405
62 3.63140781173228
63 3.79949540851117
64 3.96873369242778
65 4.14450312787326
66 4.32186529214627
67 4.4932582488657
68 4.67112433034814
69 4.84178565635781
70 5.02901545299445
71 5.21289545198924
72 5.38576933808042
73 5.46206987485165
74 5.64471375020515
75 5.81927141996669
76 5.98175574238547
77 6.14214359468317
78 6.30324701788274
79 6.47223240138462
80 6.65769850956447
81 6.84480684370334
82 7.0309720978029
83 7.22477236014694
84 7.40943009354761
85 7.6074579274753
86 7.80077098781539
87 8.00442961996107
88 8.15348379888818
89 8.28647329477318
90 8.42167744613754
91 8.55940669897098
92 8.6925647251133
93 8.81520761745009
94 8.93607201845955
95 9.10143584021288
96 9.29233486476632
97 9.46969604108904
98 9.62563892587564
99 9.69989865882668
100 9.78456964725479
101 9.96437668539115
102 10.1347431686104
103 10.3047676597742
104 10.4485219521226
105 10.6082825815518
106 10.7767072932144
107 10.9841421585569
108 11.1570753674698
109 11.3453545967615
110 11.5786543217887
111 11.7724909998503
112 11.990381482652
113 12.1886727493032
114 12.3667355225657
115 12.4682302918176
116 12.7127242235289
117 12.9433210102919
118 13.1567209149774
119 13.3837419703974
};
\addlegendentry{TOMGRO}
\end{axis}

\end{tikzpicture}
	\caption{Comparison of the \textit{SIMPLE} and \textit{TOMGRO} with the \textit{DSSAT} model.} 
	\label{fig:validation_dssat}
	\vspace{-0.2cm}
\end{figure}
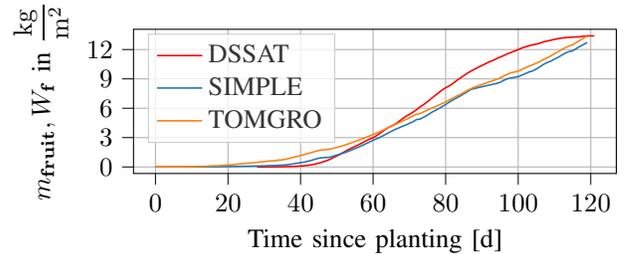

Furthermore, we verify the accuracy of both models by subjecting them to validation in the \textit{DSSAT} environment \cite{Jones.2003}, which is widely used for crop growth modeling. For simulation of tomato growth and weather data, we utilize the \textit{DSSAT} model. Fig. \ref{fig:validation_dssat} illustrates that the \textit{SIMPLE} and \textit{TOMGRO} models produce nearly identical trajectories, whereas the \textit{DSSAT} simulation behaves differently starting from day 60. However, all three models generate a similar yield of approximately \SI{13}{\kilogram\per\square\metre}.

Hence, the models demonstrate their ability to simulate real-life scenarios. 
In the following, the suitability of the models for use in \ac{OC} is investigated to find optimal inputs for the greenhouse.

\subsection{Comparing of Optimal Control Results}
In this section, we present the results of applying the \ac{OC} framework, as described in Section \ref{sec:optimalcontrol-framwork}, to the two crop models. The cost function \eqref{eq:cost} weights are chosen such that the cost represents the negative value of the economic yield, given that the fresh weight of the tomato is sold for \SI{2}{\sieuro\per\kilogram}. The cost for CO\textsubscript{2} is assumed to be \SI{0.15}{\sieuro\per\kilogram}, while ventilation is considered to be neutral in terms of cost. Heating is assumed to cost \SI{0.02}{\sieuro\per\kilo\watt\per\hour}. The results include the optimal values of the crop model states, greenhouse states, and control inputs.
Fig. \ref{fig:states_simple_co2} and Fig. \ref{fig:states_TOMGRO_co2} show the optimal biomass evolution of the \textit{SIMPLE} and \textit{TOMGRO} model, respectively. Their optimal results can be found in Tab. \ref{tab:Summary}. 
\begin{figure}[htb]
	\centering 
\begin{tikzpicture}

\definecolor{darkgray176}{RGB}{176,176,176}
\definecolor{green}{RGB}{0,128,0}
\definecolor{lightgray204}{RGB}{204,204,204}
\definecolor{magenta}{RGB}{255,0,255}

\begin{axis}[
height=3.5cm,
legend cell align={left},
legend style={
  fill opacity=0.8,
  draw opacity=1,
  text opacity=1,
  at={(0.03,0.97)},
  anchor=north west,
  draw=lightgray204
},
minor xtick={},
minor ytick={},
tick align=outside,
tick pos=left,
width=7cm,
x grid style={darkgray176},
xlabel={Time since planting in d},
xmajorgrids,
xmin=-5, xmax=105,
xtick style={color=black},
xtick={-20,0,20,40,60,80,100,120},
y grid style={darkgray176},
ylabel={\(\displaystyle {\mathrm{m_\mathrm{B}}}\) in \(\displaystyle {\mathrm{\frac{kg}{m^2}}}\)},
ymajorgrids,
ymin=-1.5, ymax=19.5,
ytick style={color=black},
ytick={0,6,12,18}
]
\addplot [semithick, green]
table {%
0 -0
1 0.00084
2 0.00195
3 0.00333
4 0.00513
5 0.00749
6 0.01056
7 0.01455
8 0.01974
9 0.02646
10 0.03514
11 0.04631
12 0.06061
13 0.07883
14 0.10186
15 0.13491
16 0.17988
17 0.236
18 0.30401
19 0.38491
20 0.47945
21 0.58798
22 0.71046
23 0.84521
24 0.98992
25 1.14385
26 1.30612
27 1.47592
28 1.65234
29 1.83459
30 2.02186
31 2.21344
32 2.40866
33 2.60696
34 2.80781
35 3.0108
36 3.21551
37 3.42164
38 3.6289
39 3.83707
40 4.04596
41 4.25542
42 4.46532
43 4.67557
44 4.88609
45 5.09682
46 5.30772
47 5.51875
48 5.72988
49 5.9411
50 6.15239
51 6.36373
52 6.57512
53 6.78654
54 6.99799
55 7.20946
56 7.42095
57 7.63245
58 7.84397
59 8.0555
60 8.26704
61 8.47857
62 8.69013
63 8.90167
64 9.11324
65 9.32479
66 9.53636
67 9.74791
68 9.95949
69 10.17104
70 10.38262
71 10.59417
72 10.80576
73 11.01731
74 11.2289
75 11.44044
76 11.65205
77 11.86357
78 12.07519
79 12.2867
80 12.49834
81 12.70983
82 12.92149
83 13.13295
84 13.34464
85 13.55605
86 13.76778
87 13.97915
88 14.19092
89 14.40221
90 14.61405
91 14.82523
92 15.03716
93 15.24819
94 15.46024
95 15.67105
96 15.88324
97 16.09373
98 16.30612
99 16.51608
100 16.72797
};
\addlegendentry{\footnotesize $m_\mathrm{B}$}
\end{axis}

\begin{axis}[
axis y line=right,
height=3.5cm,
legend cell align={left},
legend style={
  fill opacity=0.8,
  draw opacity=1,
  text opacity=1,
  at={(0.97,0.1)},
  anchor=south east,
  draw=lightgray204
},
minor xtick={},
minor ytick={},
tick align=outside,
width=7cm,
x grid style={darkgray176},
xmajorgrids,
xmin=-5, xmax=105,
xtick pos=left,
xtick style={color=black},
xtick={-20,0,20,40,60,80,100,120},
y grid style={darkgray176},
ylabel={ \(\displaystyle {\mathrm{\tau}}, {\mathrm{I_{50B}}}\) in \(\displaystyle ^\circ\)Cd},
ymajorgrids,
ymin=-200, ymax=2600,
ytick pos=right,
ytick style={color=black},
ytick={0,800,1600,2400},
yticklabel style={anchor=west}
]
\addplot [semithick, magenta]
table {%
0 -0
1 17.98958
2 43.95028
3 70.64876
4 97.40419
5 124.16406
6 150.92427
7 177.68451
8 204.44475
9 231.205
10 257.96524
11 284.72548
12 311.48572
13 338.24596
14 365.00621
15 391.76645
16 418.18118
17 444.32233
18 470.34742
19 496.27389
20 522.09387
21 547.7444
22 571.85774
23 592.73302
24 613.20937
25 633.60633
26 653.90975
27 674.18019
28 694.40623
29 714.61455
30 734.7975
31 754.97009
32 775.12703
33 795.27776
34 815.41851
35 835.55555
36 855.68616
37 875.81462
38 895.93895
39 916.0621
40 936.18263
41 956.30253
42 976.42079
43 996.53871
44 1016.65562
45 1036.7723
46 1056.8884
47 1077.00429
48 1097.11992
49 1117.23527
50 1137.35063
51 1157.4656
52 1177.58082
53 1197.6955
54 1217.81069
55 1237.92512
56 1258.04035
57 1278.15454
58 1298.26989
59 1318.38381
60 1338.49935
61 1358.61297
62 1378.7288
63 1398.84203
64 1418.95826
65 1439.07099
66 1459.18777
67 1479.29984
68 1499.41736
69 1519.52857
70 1539.6471
71 1559.75714
72 1579.87702
73 1599.98551
74 1620.10722
75 1640.21363
76 1660.33779
77 1680.44143
78 1700.56891
79 1720.66884
80 1740.80081
81 1760.89578
82 1781.03382
83 1801.1222
84 1821.2685
85 1841.3481
86 1861.50567
87 1881.57362
88 1901.74674
89 1921.79923
90 1941.99404
91 1962.02602
92 1982.25177
93 2002.25648
94 2022.52797
95 2042.49602
96 2062.83985
97 2082.75745
98 2103.2329
99 2123.07877
100 2143.54161
};
\addlegendentry{\footnotesize ${\mathrm{\tau}}$}
\addplot [semithick, black]
table {%
0 -0
1 0.00691
2 0.05713
3 0.33672
4 0.63738
5 0.9397
6 1.24214
7 1.54459
8 1.84705
9 2.1495
10 2.45196
11 2.75441
12 3.05686
13 3.35932
14 3.66177
15 3.96422
16 4.14139
17 4.23126
18 4.29325
19 4.33832
20 4.37157
21 4.39469
22 4.40411
23 4.41144
24 4.41869
25 4.42592
26 4.43313
27 4.44034
28 4.44754
29 4.45473
30 4.46193
31 4.46911
32 4.4763
33 4.48349
34 4.49067
35 4.49785
36 4.50503
37 4.51221
38 4.51939
39 4.52657
40 4.53375
41 4.54093
42 4.54811
43 4.55529
44 4.56247
45 4.56965
46 4.57683
47 4.58401
48 4.59118
49 4.59836
50 4.60554
51 4.61272
52 4.6199
53 4.62708
54 4.63426
55 4.64143
56 4.64861
57 4.65579
58 4.66297
59 4.67015
60 4.67733
61 4.68451
62 4.69168
63 4.69886
64 4.70604
65 4.71322
66 4.7204
67 4.72758
68 4.73476
69 4.74193
70 4.74911
71 4.75629
72 4.76347
73 4.77065
74 4.77783
75 4.785
76 4.79218
77 4.79936
78 4.80654
79 4.81372
80 4.8209
81 4.82807
82 4.83526
83 4.84243
84 4.84961
85 4.85679
86 4.86397
87 4.87114
88 4.87833
89 4.8855
90 4.89269
91 4.89986
92 4.90705
93 4.91421
94 4.92142
95 4.92857
96 4.9358
97 4.94294
98 4.95019
99 4.95732
100 4.96456
};
\addlegendentry{\footnotesize ${\mathrm{I_{50B}}}$}
\end{axis}

\end{tikzpicture}
	\vspace{-0.5cm}
	\caption{Optimized states of the \textit{SIMPLE}  model.} 
	\label{fig:states_simple_co2}
	\vspace{-0.2cm}
\end{figure}
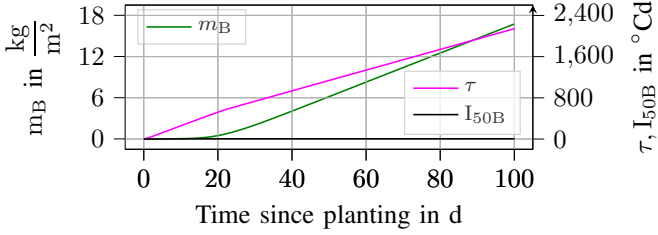  
\begin{figure}[htb]
	\centering 
\begin{tikzpicture}

\definecolor{darkgray176}{RGB}{176,176,176}
\definecolor{green}{RGB}{0,128,0}
\definecolor{lightgray204}{RGB}{204,204,204}
\definecolor{magenta}{RGB}{255,0,255}

\begin{groupplot}[group style={group size=1 by 2}]
\nextgroupplot[
height=3.5cm,
legend cell align={left},
legend style={
  fill opacity=0.8,
  draw opacity=1,
  text opacity=1,
  at={(0.03,0.97)},
  anchor=north west,
  draw=lightgray204
},
minor xtick={},
minor ytick={},
scaled x ticks=manual:{}{\pgfmathparse{#1}},
tick align=outside,
tick pos=left,
width=7cm,
x grid style={darkgray176},
xmajorgrids,
xmin=-5, xmax=105,
xtick style={color=black},
xtick={-20,0,20,40,60,80,100,120},
xticklabels={},
y grid style={darkgray176},
ylabel={\(\displaystyle N\)},
ymajorgrids,
ymin=-6.25, ymax=56.25,
ytick style={color=black},
ytick={0,25,50}
]
\addplot [semithick, green]
table {%
0 6
1 6.46002
2 6.95837
3 7.4561
4 7.95439
5 8.45214
6 8.9503
7 9.44802
8 9.94601
9 10.44364
10 10.94144
11 11.43893
12 11.93648
13 12.43374
14 12.93096
15 13.4279
16 13.92472
17 14.42121
18 14.91751
19 15.41342
20 15.90904
21 16.40418
22 16.89891
23 17.39301
24 17.88655
25 18.37924
26 18.87115
27 19.36189
28 19.85152
29 20.33939
30 20.82991
31 21.32434
32 21.82011
33 22.31614
34 22.81235
35 23.308
36 23.8033
37 24.29735
38 24.79064
39 25.28179
40 25.77177
41 26.25808
42 26.74059
43 27.20748
44 27.66546
45 28.12319
46 28.57909
47 29.03535
48 29.49064
49 29.94608
50 30.40094
51 30.85581
52 31.31029
53 31.7647
54 32.21883
55 32.67285
56 33.12664
57 33.58029
58 34.03373
59 34.48703
60 34.94013
61 35.39307
62 35.84582
63 36.29839
64 36.75077
65 37.20297
66 37.65496
67 38.10675
68 38.55829
69 39.00963
70 39.46066
71 39.9115
72 40.36192
73 40.8122
74 41.26183
75 41.71147
76 42.15996
77 42.60888
78 43.05536
79 43.50338
80 43.94515
81 44.3915
82 44.80262
83 45.23946
84 45.6745
85 46.10855
86 46.53763
87 46.89662
88 47.29568
89 47.7005
90 48.05565
91 48.4011
92 48.74502
93 49.0887
94 49.43233
95 49.77595
96 50.11957
97 50.46319
98 50.80681
99 51.15043
100 51.5899
};
\addlegendentry{\footnotesize $N$}

\nextgroupplot[
height=3.5cm,
legend cell align={left},
legend style={
  fill opacity=0.8,
  draw opacity=1,
  text opacity=1,
  at={(0.03,0.97)},
  anchor=north west,
  draw=lightgray204
},
minor xtick={},
minor ytick={},
tick align=outside,
tick pos=left,
width=7cm,
x grid style={darkgray176},
xlabel={Time since planting in d},
xmajorgrids,
xmin=-5, xmax=105,
xtick style={color=black},
xtick={-20,0,20,40,60,80,100,120},
y grid style={darkgray176},
ylabel={\(\displaystyle W_f,W\) in \(\displaystyle \mathrm{\frac{kg}{m^2}}\)},
ymajorgrids,
ymin=-2.5, ymax=22.5,
ytick style={color=black},
ytick={0,10,20}
]
\addplot [semithick, blue]
table {%
0 0
1 -0
2 -0
3 -1e-05
4 -2e-05
5 -3e-05
6 -4e-05
7 -7e-05
8 -0.0001
9 -0.00015
10 -0.00021
11 -0.0003
12 -0.00042
13 -0.00059
14 -0.00082
15 -0.00111
16 -0.00151
17 -0.00202
18 -0.00269
19 -0.00353
20 -0.00459
21 -0.00591
22 -0.00754
23 -0.00952
24 -0.01188
25 -0.01465
26 -0.01781
27 -0.02132
28 -0.02507
29 -0.02894
30 -0.03234
31 -0.03489
32 -0.0365
33 -0.03698
34 -0.03606
35 -0.03345
36 -0.02891
37 -0.02202
38 -0.01274
39 -0.00037
40 0.01475
41 0.0341
42 0.05888
43 0.09674
44 0.14294
45 0.1954
46 0.25361
47 0.31734
48 0.38642
49 0.46065
50 0.53987
51 0.62387
52 0.7125
53 0.80557
54 0.9029
55 1.00433
56 1.10968
57 1.2188
58 1.33152
59 1.4477
60 1.56717
61 1.68979
62 1.81534
63 1.94349
64 2.07412
65 2.20711
66 2.34236
67 2.47973
68 2.61915
69 2.76048
70 2.90366
71 3.04855
72 3.19511
73 3.34319
74 3.4928
75 3.64373
76 3.7961
77 3.94956
78 4.1045
79 4.26019
80 4.41775
81 4.57549
82 4.73831
83 4.89878
84 5.06028
85 5.22248
86 5.38597
87 5.55761
88 5.72551
89 5.89316
90 6.06651
91 6.24103
92 6.41594
93 6.59108
94 6.76643
95 6.94196
96 7.11764
97 7.29344
98 7.46935
99 7.64532
100 7.80373
};
\addlegendentry{\footnotesize $W_f$}
\addplot [semithick, black]
table {%
0 0.004
1 0.00047
2 -0.00158
3 -0.00144
4 0.00103
5 0.00597
6 0.01352
7 0.02386
8 0.03711
9 0.05348
10 0.07309
11 0.09695
12 0.12602
13 0.16066
14 0.20116
15 0.24787
16 0.30107
17 0.36108
18 0.42818
19 0.50266
20 0.58479
21 0.67481
22 0.77297
23 0.87948
24 0.99453
25 1.11833
26 1.25036
27 1.39023
28 1.53799
29 1.69377
30 1.85689
31 2.02227
32 2.18898
33 2.35681
34 2.52598
35 2.69654
36 2.86878
37 3.04281
38 3.21884
39 3.39709
40 3.57756
41 3.76086
42 3.94818
43 4.14315
44 4.34331
45 4.54948
46 4.76064
47 4.97728
48 5.19877
49 5.4253
50 5.65645
51 5.89221
52 6.13229
53 6.37661
54 6.62491
55 6.87708
56 7.13291
57 7.39226
58 7.65496
59 7.92085
60 8.18978
61 8.4616
62 8.73607
63 9.01287
64 9.29188
65 9.57298
66 9.85604
67 10.14094
68 10.42758
69 10.71584
70 11.00562
71 11.29682
72 11.58933
73 11.88307
74 12.17789
75 12.4738
76 12.77054
77 13.06826
78 13.36644
79 13.66563
80 13.96443
81 14.26465
82 14.55847
83 14.8579
84 15.15757
85 15.45738
86 15.75667
87 16.04173
88 16.33553
89 16.63069
90 16.91573
91 17.19872
92 17.48147
93 17.76424
94 18.04704
95 18.32986
96 18.61269
97 18.89549
98 19.17825
99 19.46093
100 19.75565
};
\addlegendentry{\footnotesize $W$}
\end{groupplot}

\begin{groupplot}[group style={group size=1 by 2}]
\nextgroupplot[
axis y line=right,
height=3.5cm,
legend cell align={left},
legend style={
  fill opacity=0.8,
  draw opacity=1,
  text opacity=1,
  at={(0.97,0.03)},
  anchor=south east,
  draw=lightgray204
},
minor xtick={},
minor ytick={},
tick align=outside,
width=7cm,
x grid style={darkgray176},
xmin=-5, xmax=105,
xtick pos=left,
xtick style={color=black},
xtick={-20,0,20,40,60,80,100,120},
y grid style={darkgray176},
ylabel={\(\displaystyle LAI\)},
ymin=-0.625, ymax=5.625,
ytick pos=right,
ytick style={color=black},
ytick={0,2.5,5},
yticklabel style={anchor=west}
]
\addplot [semithick, red]
table {%
0 0.006
1 0.01441
2 0.02414
3 0.03456
4 0.04574
5 0.05768
6 0.07046
7 0.08409
8 0.09863
9 0.1141
10 0.13055
11 0.14801
12 0.16651
13 0.18609
14 0.20678
15 0.22859
16 0.25155
17 0.27569
18 0.30101
19 0.32752
20 0.35523
21 0.38414
22 0.41424
23 0.44551
24 0.47795
25 0.51152
26 0.54621
27 0.58195
28 0.61873
29 0.65647
30 0.69546
31 0.73582
32 0.77729
33 0.81978
34 0.86324
35 0.90756
36 0.95272
37 0.99859
38 1.04519
39 1.09232
40 1.14005
41 1.18807
42 1.23633
43 1.28359
44 1.33043
45 1.37772
46 1.42524
47 1.47321
48 1.52146
49 1.57008
50 1.61897
51 1.66818
52 1.71764
53 1.76737
54 1.81733
55 1.86751
56 1.91789
57 1.96847
58 2.01922
59 2.07013
60 2.12119
61 2.17239
62 2.22372
63 2.27516
64 2.3267
65 2.37834
66 2.43007
67 2.48187
68 2.53374
69 2.58567
70 2.63764
71 2.68967
72 2.74171
73 2.7938
74 2.84587
75 2.89799
76 2.95003
77 3.00216
78 3.05404
79 3.10614
80 3.15754
81 3.2095
82 3.25736
83 3.30825
84 3.35894
85 3.40952
86 3.45952
87 3.50132
88 3.5478
89 3.59494
90 3.63624
91 3.67637
92 3.71627
93 3.75608
94 3.7958
95 3.83538
96 3.87473
97 3.91363
98 3.95142
99 3.98398
100 3.99678
};
\addlegendentry{\footnotesize $LAI$}

\nextgroupplot[
axis y line=right,
height=3.5cm,
legend cell align={left},
legend style={
  fill opacity=0.8,
  draw opacity=1,
  text opacity=1,
  at={(0.03,0.2)},
  anchor=south west,
  draw=lightgray204
},
minor xtick={},
minor ytick={},
tick align=outside,
width=7cm,
x grid style={darkgray176},
xmin=-5, xmax=105,
xtick pos=left,
xtick style={color=black},
xtick={-20,0,20,40,60,80,100,120},
y grid style={darkgray176},
ylabel={\(\displaystyle W_m\) in \(\displaystyle \mathrm{\frac{kg}{m^2}}\)},
ymin=-0.375, ymax=3.375,
ytick pos=right,
ytick style={color=black},
ytick={0,1.5,3},
yticklabel style={anchor=west}
]
\addplot [semithick, magenta]
table {%
0 0
1 0
2 -0
3 -0
4 -0
5 -0
6 -0
7 -0
8 -0
9 -0
10 -0
11 -0
12 -0
13 -0
14 -0
15 -0
16 -0
17 -0
18 -0
19 -0
20 -0
21 -0
22 -0
23 -0
24 -0
25 -1e-05
26 -1e-05
27 -2e-05
28 -3e-05
29 -5e-05
30 -7e-05
31 -0.00011
32 -0.00015
33 -0.00021
34 -0.00029
35 -0.00039
36 -0.0005
37 -0.00062
38 -0.00072
39 -0.0008
40 -0.00079
41 -0.00066
42 -0.00032
43 0.00031
44 0.00146
45 0.00332
46 0.00607
47 0.0099
48 0.01497
49 0.02146
50 0.02948
51 0.03916
52 0.05061
53 0.06392
54 0.07913
55 0.09633
56 0.11554
57 0.1368
58 0.16015
59 0.18561
60 0.21318
61 0.24289
62 0.27475
63 0.30876
64 0.34492
65 0.38324
66 0.42369
67 0.46628
68 0.51099
69 0.55783
70 0.60676
71 0.6578
72 0.71089
73 0.76606
74 0.82323
75 0.88246
76 0.94359
77 1.00681
78 1.0717
79 1.1388
80 1.20695
81 1.27775
82 1.34487
83 1.41818
84 1.49308
85 1.56968
86 1.64724
87 1.71386
88 1.7898
89 1.86868
90 1.9397
91 2.01078
92 2.08338
93 2.15768
94 2.23371
95 2.31144
96 2.39085
97 2.47191
98 2.55459
99 2.63887
100 2.74767
};
\addlegendentry{\footnotesize $W_m$}
\end{groupplot}

\end{tikzpicture}
	\vspace{-0.5cm}
	\caption{Optimized states of the \textit{TOMGRO} model.} 
	\vspace{-0.2cm}
	\label{fig:states_TOMGRO_co2}
	\vspace{-0.3cm}
\end{figure}
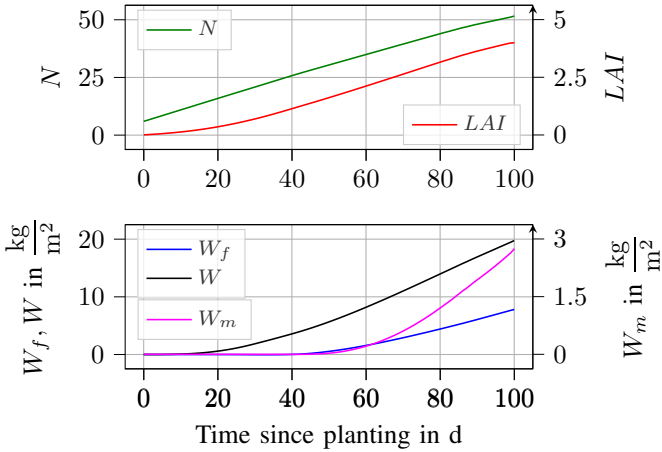  
Fig. \ref{fig:greenhouse_inputs_simple_co2} displays the \ac{OC} inputs of the greenhouse. The optimal value for $u^{Ap}_{v,i}$ is 0 throughout for both crop models since ventilation leads to a decrease in the CO\textsubscript{2} concentration. 
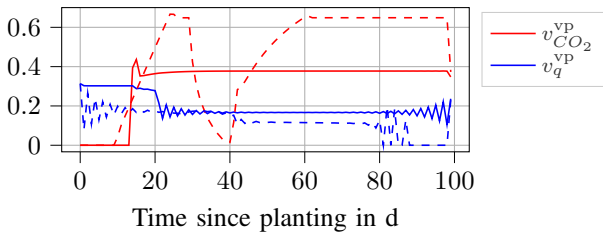
\begin{figure}[tb]
	\centering 
	\vspace{-0.5cm}
\begin{tikzpicture}

\definecolor{darkgray176}{RGB}{176,176,176}
\definecolor{lightgray204}{RGB}{204,204,204}

\begin{axis}[
height=3.5cm,
legend cell align={left},
legend style={
  fill opacity=0.8,
  draw opacity=1,
  text opacity=1,
  at={(1.03,0.97)},
  anchor=north west,
  draw=lightgray204
},
minor xtick={},
minor ytick={},
tick align=outside,
tick pos=left,
width=7cm,
x grid style={darkgray176},
xlabel={Time since planting in d},
xmajorgrids,
xmin=-4.95, xmax=103.95,
xtick style={color=black},
xtick={-20,0,20,40,60,80,100,120},
y grid style={darkgray176},
ymajorgrids,
ymin=-0.0333495, ymax=0.7003395,
ytick style={color=black},
ytick={0,0.2,0.4,0.6}
]
\addplot [semithick, red]
table {%
0 0
1 0
2 0
3 0
4 0
5 0
6 0
7 0
8 0
9 0
10 0
11 0
12 0
13 0
14 0.39394
15 0.43602
16 0.35245
17 0.35443
18 0.35797
19 0.36127
20 0.36418
21 0.3667
22 0.36845
23 0.36966
24 0.37097
25 0.37207
26 0.37303
27 0.37382
28 0.3745
29 0.37506
30 0.37554
31 0.37593
32 0.37625
33 0.37651
34 0.37673
35 0.3769
36 0.37704
37 0.37714
38 0.37723
39 0.3773
40 0.37735
41 0.37739
42 0.37742
43 0.37744
44 0.37746
45 0.37748
46 0.37749
47 0.3775
48 0.37751
49 0.37752
50 0.37752
51 0.37753
52 0.37753
53 0.37753
54 0.37753
55 0.37754
56 0.37754
57 0.37754
58 0.37754
59 0.37754
60 0.37754
61 0.37754
62 0.37754
63 0.37754
64 0.37754
65 0.37754
66 0.37754
67 0.37755
68 0.37754
69 0.37755
70 0.37754
71 0.37755
72 0.37753
73 0.37755
74 0.37753
75 0.37755
76 0.37753
77 0.37756
78 0.37752
79 0.37756
80 0.37752
81 0.37757
82 0.37751
83 0.37758
84 0.37749
85 0.3776
86 0.37748
87 0.37761
88 0.37746
89 0.37764
90 0.37742
91 0.37767
92 0.37738
93 0.37771
94 0.37732
95 0.37777
96 0.37723
97 0.37785
98 0.37709
99 0.34769
};
\addlegendentry{\footnotesize $v^\mathrm{vp}_{CO_2}$}
\addplot [semithick, blue]
table {%
0 0.31355
1 0.30288
2 0.3024
3 0.30236
4 0.30236
5 0.30236
6 0.30236
7 0.30236
8 0.30236
9 0.30236
10 0.30236
11 0.30236
12 0.30236
13 0.30236
14 0.30236
15 0.28786
16 0.28996
17 0.28467
18 0.28541
19 0.2808
20 0.27798
21 0.21595
22 0.1363
23 0.2052
24 0.14375
25 0.19336
26 0.14936
27 0.1858
28 0.15361
29 0.18056
30 0.15682
31 0.17683
32 0.15925
33 0.17413
34 0.16108
35 0.17217
36 0.16247
37 0.17073
38 0.16352
39 0.16968
40 0.16431
41 0.16891
42 0.1649
43 0.16835
44 0.16535
45 0.16794
46 0.16568
47 0.16763
48 0.16592
49 0.16741
50 0.1661
51 0.16725
52 0.16622
53 0.16714
54 0.16631
55 0.16707
56 0.16636
57 0.16704
58 0.16638
59 0.16703
60 0.16637
61 0.16705
62 0.16633
63 0.1671
64 0.16626
65 0.16719
66 0.16616
67 0.16732
68 0.166
69 0.1675
70 0.16578
71 0.16776
72 0.16549
73 0.1681
74 0.16509
75 0.16857
76 0.16454
77 0.16921
78 0.16381
79 0.17006
80 0.16283
81 0.17121
82 0.16151
83 0.17275
84 0.15973
85 0.17484
86 0.15734
87 0.17766
88 0.15412
89 0.18149
90 0.14977
91 0.18674
92 0.14382
93 0.19405
94 0.13559
95 0.20455
96 0.12379
97 0.22078
98 0.10569
99 0.23646
};
\addlegendentry{\footnotesize $v^\mathrm{vp}_q$}
\addplot [semithick, red, dashed, forget plot]
table {%
0 0
1 0
2 0
3 0
4 0
5 0
6 0
7 0
8 0
9 0
10 0.04333
11 0.10525
12 0.14888
13 0.19219
14 0.23588
15 0.2793
16 0.323
17 0.36644
18 0.41003
19 0.45331
20 0.49662
21 0.53954
22 0.5824
23 0.62477
24 0.66699
25 0.66621
26 0.64858
27 0.64858
28 0.64859
29 0.6486
30 0.42416
31 0.30901
32 0.23462
33 0.18481
34 0.13825
35 0.10644
36 0.07472
37 0.05303
38 0.03098
39 0.01597
40 0.01431
41 0.11536
42 0.28288
43 0.29809
44 0.32409
45 0.35498
46 0.38559
47 0.41381
48 0.44015
49 0.46468
50 0.48756
51 0.50895
52 0.52898
53 0.54776
54 0.56541
55 0.58201
56 0.59765
57 0.61241
58 0.62635
59 0.63954
60 0.65202
61 0.65643
62 0.64878
63 0.64878
64 0.64878
65 0.64878
66 0.64878
67 0.64879
68 0.64879
69 0.64879
70 0.64879
71 0.64879
72 0.6488
73 0.6488
74 0.6488
75 0.6488
76 0.6488
77 0.6488
78 0.6488
79 0.6488
80 0.64881
81 0.64881
82 0.64882
83 0.64881
84 0.64881
85 0.64881
86 0.64881
87 0.64883
88 0.64882
89 0.64882
90 0.64883
91 0.64883
92 0.64883
93 0.64883
94 0.64883
95 0.64883
96 0.64883
97 0.64883
98 0.64883
99 0.37307
};
\addplot [semithick, blue, dashed, forget plot]
table {%
0 0.31355
1 0.07944
2 0.26364
3 0.10861
4 0.23858
5 0.12927
6 0.22051
7 0.14375
8 0.20758
9 0.15372
10 0.19823
11 0.16035
12 0.19133
13 0.16459
14 0.18605
15 0.16709
16 0.18182
17 0.16826
18 0.17818
19 0.16838
20 0.17477
21 0.16761
22 0.1713
23 0.166
24 0.16744
25 0.16355
26 0.16287
27 0.16005
28 0.15691
29 0.1697
30 0.17301
31 0.17436
32 0.17381
33 0.175
34 0.1711
35 0.17316
36 0.1659
37 0.16974
38 0.1583
39 0.16512
40 0.14719
41 0.15278
42 0.10659
43 0.13127
44 0.11126
45 0.12476
46 0.11408
47 0.12119
48 0.11547
49 0.11915
50 0.11604
51 0.1179
52 0.11617
53 0.11706
54 0.11606
55 0.11645
56 0.11583
57 0.11594
58 0.11552
59 0.11549
60 0.11517
61 0.11505
62 0.1148
63 0.1146
64 0.1144
65 0.11412
66 0.11396
67 0.11356
68 0.11352
69 0.11287
70 0.11308
71 0.11193
72 0.11276
73 0.11046
74 0.1128
75 0.10784
76 0.1138
77 0.10228
78 0.11731
79 0.08781
80 0.12799
81 0
82 0.19764
83 0
84 0.18248
85 0
86 0
87 0.12872
88 0
89 0
90 0
91 0
92 0
93 0
94 0
95 0
96 0
97 0
98 -0
99 0.22929
};
\end{axis}

\end{tikzpicture}
	\caption{Optimal inputs of the \ac{GH} with the use of different tomato models in the optimization: Solid lines: \textit{SIMPLE}  model; Dashed lines: \textit{TOMGRO} model.}
	\label{fig:greenhouse_inputs_simple_co2}
	\vspace{-0.0cm}
\end{figure}  
Therefore, for simplicity, ventilation is not additionally marked in the graph. We conclude that in the optimal solution, the greenhouse ventilation is not used, indicating that the optimization aims to maintain a stable environment throughout the cultivation period. The optimal values of $u^{vp}_{q}$ and $u^{vp}_{CO_2}$ of the greenhouse show similar results.

\begin{figure}[ht]
    \vspace{-0.0cm}
	\centering 
	\input{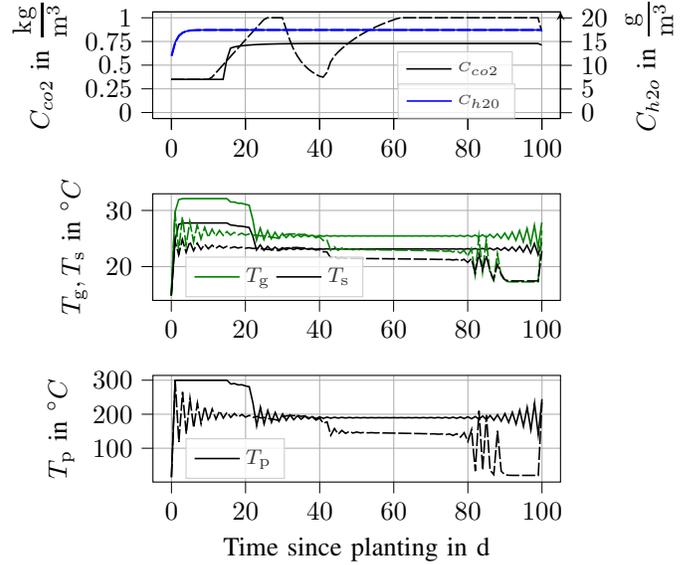}
	\vspace{-0.5cm}
	\caption{Optimal sequence of states of the \ac{GH}: Solid lines: \textit{SIMPLE}  model; Dashed lines: \textit{TOMGRO} model.} 
	\label{fig:greenhouse_states_simple_co2}
	\vspace{-0.2cm}
\end{figure}

The optimal states of the greenhouse are shown in Fig. \ref{fig:greenhouse_states_simple_co2}. We obtain similar results for both crop models, particularly for the $T_s$, $T_p$, and $C_{H_{2}O}$. The optimized temperature $T_g$ for the \textit{SIMPLE} model is around \SI{27}{\degreeCelsius} during the cultivation period, which is close to the optimal temperature for this tomato variety (\SI{26}{\degreeCelsius}) according to \cite{Zhao.2019}. The optimal solution of $C_{CO_{2}}$ concentration differs between the models, with the concentration resulting from the \textit{SIMPLE} model saturating at \SI{700}{ppm}, since further increase makes no difference in the simplified model.

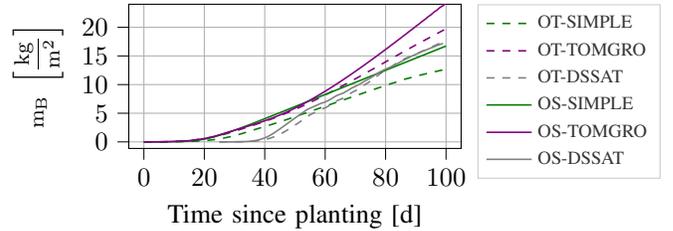
\begin{figure}[ht]
	\centering 
\begin{tikzpicture}

\definecolor{darkgray176}{RGB}{176,176,176}
\definecolor{gray}{RGB}{128,128,128}
\definecolor{green}{RGB}{0,128,0}
\definecolor{lightgray204}{RGB}{204,204,204}
\definecolor{purple}{RGB}{128,0,128}

\begin{axis}[
height=3.5cm,
legend cell align={left},
legend style={
  fill opacity=0.8,
  draw opacity=1,
  text opacity=1,
  at={(1.6, 1.0)},
  draw=lightgray204
},
minor xtick={},
minor ytick={},
tick align=outside,
tick pos=left,
unbounded coords=jump,
width=6cm,
x grid style={darkgray176},
xlabel={Time since planting [d]},
xmajorgrids,
xmin=-5, xmax=105,
xtick style={color=black},
xtick={-20,0,20,40,60,80,100,120},
y grid style={darkgray176},
ylabel={\footnotesize \(\displaystyle {\mathrm{m_\mathrm{B}}}\ \left[{\mathrm{\frac{kg}{m^2}}}\right]\)},
ymajorgrids,
ymin=-1.1455175, ymax=24.0206675,
ytick style={color=black},
ytick={-5,0,5,10,15,20,25}
]
\addplot [semithick, green, dashed]
table {%
0 0
1 0.00042
2 0.00144
3 0.00254
4 0.00408
5 0.00582
6 0.00815
7 0.01083
8 0.01431
9 0.0184
10 0.02358
11 0.02972
12 0.0375
13 0.04704
14 0.05911
15 0.07389
16 0.09237
17 0.11485
18 0.14259
19 0.17593
20 0.21605
21 0.26281
22 0.31756
23 0.38032
24 0.4524
25 0.5336
26 0.62489
27 0.72581
28 0.83666
29 0.95688
30 1.0858
31 1.22695
32 1.37824
33 1.53811
34 1.70242
35 1.86704
36 2.03078
37 2.19493
38 2.35733
39 2.52047
40 2.68037
41 2.84168
42 2.99829
43 3.16234
44 3.32286
45 3.49637
46 3.6685
47 3.84865
48 4.02807
49 4.21086
50 4.39234
51 4.57535
52 4.75759
53 4.94061
54 5.12321
55 5.30621
56 5.48895
57 5.67188
58 5.85465
59 6.03748
60 6.22019
61 6.40289
62 6.58549
63 6.76804
64 6.95049
65 7.13286
66 7.31513
67 7.49728
68 7.67935
69 7.86124
70 8.04308
71 8.22465
72 8.40625
73 8.58738
74 8.76878
75 8.94926
76 9.13054
77 9.3099
78 9.49138
79 9.66844
80 9.85096
81 10.02214
82 10.20799
83 10.34433
84 10.55173
85 10.69303
86 10.89477
87 11.03489
88 11.15988
89 11.33616
90 11.4715
91 11.59588
92 11.7186
93 11.84105
94 11.96346
95 12.08586
96 12.20826
97 12.33066
98 12.45306
99 12.57546
};
\addlegendentry{\scriptsize OT-SIMPLE}
\addplot [semithick, purple, dashed]
table {%
0 0.004
1 0.00047
2 -0.00158
3 -0.00144
4 0.00103
5 0.00597
6 0.01352
7 0.02386
8 0.03711
9 0.05348
10 0.07309
11 0.09695
12 0.12602
13 0.16066
14 0.20116
15 0.24787
16 0.30107
17 0.36108
18 0.42818
19 0.50266
20 0.58479
21 0.67481
22 0.77297
23 0.87948
24 0.99453
25 1.11833
26 1.25036
27 1.39023
28 1.53799
29 1.69377
30 1.85689
31 2.02227
32 2.18898
33 2.35681
34 2.52598
35 2.69654
36 2.86878
37 3.04281
38 3.21884
39 3.39709
40 3.57756
41 3.76086
42 3.94818
43 4.14315
44 4.34331
45 4.54948
46 4.76064
47 4.97728
48 5.19877
49 5.4253
50 5.65645
51 5.89221
52 6.13229
53 6.37661
54 6.62491
55 6.87708
56 7.13291
57 7.39226
58 7.65496
59 7.92085
60 8.18978
61 8.4616
62 8.73607
63 9.01287
64 9.29189
65 9.57298
66 9.85604
67 10.14094
68 10.42758
69 10.71585
70 11.00562
71 11.29683
72 11.58933
73 11.88307
74 12.17789
75 12.4738
76 12.77054
77 13.06826
78 13.36644
79 13.66563
80 13.96444
81 14.26465
82 14.55847
83 14.8579
84 15.15757
85 15.45738
86 15.75667
87 16.04173
88 16.33553
89 16.63069
90 16.91573
91 17.19873
92 17.48148
93 17.76424
94 18.04704
95 18.32986
96 18.61269
97 18.89549
98 19.17825
99 19.46093
100 19.75565
};
\addlegendentry{\scriptsize OT-TOMGRO}
\addplot [semithick, gray, dashed]
table {%
25 0
26 0
27 0.0005
28 0.0018
29 0.0042
30 0.0078
31 0.0131
32 0.0205
33 0.0304
34 0.0432
35 0.0676
36 0.1066
37 0.1478
38 0.2095
39 0.2882
40 0.3958
41 0.5253
42 0.6775
43 0.8808
44 1.1102
45 1.3392
46 1.605
47 1.9137
48 2.242
49 2.5776
50 2.9198
51 3.2746
52 3.6248
53 3.9764
54 4.3427
55 4.74
56 4.9867
57 5.252
58 5.5208
59 5.7381
60 5.9407
61 6.2132
62 6.5466
63 6.8925
64 7.2015
65 7.4188
66 7.6321
67 7.9698
68 8.317
69 8.668
70 8.9212
71 9.22
72 9.5552
73 9.945
74 10.265
75 10.6228
76 11.0195
77 11.4108
78 11.8199
79 12.2096
80 12.5552
81 12.8038
82 13.128
83 13.5172
84 13.8747
85 14.1746
86 14.4477
87 14.6704
88 14.8673
89 15.0754
90 15.3025
91 15.5421
92 15.8053
93 16.0153
94 16.2613
95 16.5277
96 16.7463
97 17.0035
98 17.2153
99 17.4263
};
\addlegendentry{\scriptsize OT-DSSAT}
\addplot [semithick, green]
table {%
0 0
1 0.00084
2 0.00195
3 0.00333
4 0.00513
5 0.00749
6 0.01056
7 0.01455
8 0.01974
9 0.02646
10 0.03514
11 0.04631
12 0.06061
13 0.07883
14 0.10186
15 0.13491
16 0.17988
17 0.236
18 0.30401
19 0.38491
20 0.47945
21 0.58798
22 0.71046
23 0.84521
24 0.98992
25 1.14385
26 1.30612
27 1.47592
28 1.65234
29 1.83459
30 2.02186
31 2.21344
32 2.40866
33 2.60696
34 2.80781
35 3.0108
36 3.21551
37 3.42164
38 3.6289
39 3.83707
40 4.04596
41 4.25542
42 4.46532
43 4.67557
44 4.88609
45 5.09682
46 5.30772
47 5.51875
48 5.72988
49 5.9411
50 6.15239
51 6.36373
52 6.57512
53 6.78654
54 6.99799
55 7.20946
56 7.42095
57 7.63245
58 7.84397
59 8.0555
60 8.26704
61 8.47857
62 8.69013
63 8.90167
64 9.11324
65 9.32479
66 9.53636
67 9.74791
68 9.95949
69 10.17104
70 10.38262
71 10.59417
72 10.80576
73 11.01731
74 11.2289
75 11.44044
76 11.65205
77 11.86357
78 12.07519
79 12.2867
80 12.49834
81 12.70983
82 12.92149
83 13.13295
84 13.34464
85 13.55605
86 13.76778
87 13.97915
88 14.19092
89 14.40221
90 14.61405
91 14.82523
92 15.03716
93 15.2482
94 15.46024
95 15.67105
96 15.88324
97 16.09373
98 16.30612
99 16.51608
100 16.72797
};
\addlegendentry{\scriptsize OS-SIMPLE}
\addplot [semithick, purple]
table {%
0 0.004
1 0.00319
2 0.00137
3 0.00211
4 0.00505
5 0.01031
6 0.018
7 0.02824
8 0.04116
9 0.05687
10 0.0755
11 0.09719
12 0.12206
13 0.15027
14 0.18195
15 0.21724
16 0.27
17 0.33802
18 0.41303
19 0.49645
20 0.58721
21 0.68695
22 0.7958
23 0.92654
24 1.07482
25 1.23067
26 1.38895
27 1.55477
28 1.71225
29 1.87559
30 2.03225
31 2.19355
32 2.34997
33 2.51042
34 2.66816
35 2.83
36 2.99171
37 3.15819
38 3.32734
39 3.50221
40 3.6823
41 3.869
42 4.06284
43 4.26386
44 4.47314
45 4.68982
46 4.91523
47 5.14796
48 5.38941
49 5.63796
50 5.89493
51 6.15866
52 6.43039
53 6.70851
54 6.99412
55 7.28573
56 7.58435
57 7.88853
58 8.19925
59 8.5151
60 8.83707
61 9.16371
62 9.49609
63 9.83265
64 10.17467
65 10.52032
66 10.87123
67 11.22516
68 11.58427
69 11.94569
70 12.31237
71 12.6805
72 13.05416
73 13.42825
74 13.80839
75 14.18766
76 14.57387
77 14.95754
78 15.34948
79 15.73673
80 16.13421
81 16.52416
82 16.9271
83 17.31878
84 17.72732
85 18.11957
86 18.53407
87 18.92554
88 19.34666
89 19.73566
90 20.16448
91 20.54882
92 20.98694
93 21.36369
94 21.81345
95 22.17844
96 22.64271
97 22.98929
98 23.46033
99 23.77597
100 24.2451
};
\addlegendentry{\scriptsize OS-TOMGRO}
\addplot [semithick, gray]
table {%
25 0
26 0
27 0.0008
28 0.0029
29 0.0065
30 0.0122
31 0.0203
32 0.0316
33 0.0464
34 0.0652
35 0.1016
36 0.1598
37 0.2433
38 0.358
39 0.5073
40 0.6987
41 0.9317
42 1.2142
43 1.5407
44 1.8959
45 2.2622
46 2.6352
47 3.0137
48 3.3935
49 3.7743
50 4.1525
51 4.5293
52 4.9001
53 5.2676
54 5.6281
55 5.883
56 6.1029
57 6.3476
58 6.5913
59 6.7866
60 6.9618
61 7.2097
62 7.5171
63 7.838
64 8.1247
65 8.328
66 8.5182
67 8.8248
68 9.1392
69 9.4582
70 9.6861
71 9.9581
72 10.2512
73 10.5919
74 10.8597
75 11.1628
76 11.5045
77 11.8414
78 12.1654
79 12.4449
80 12.7041
81 12.9124
82 13.1493
83 13.4106
84 13.6906
85 13.9702
86 14.2621
87 14.5321
88 14.7892
89 15.0589
90 15.2794
91 15.5287
92 15.7753
93 15.958
94 16.1788
95 16.4008
96 16.595
97 16.8233
98 17.0117
99 17.1914
};
\addlegendentry{\scriptsize OS-DSSAT}
\end{axis}

\end{tikzpicture}
  	\caption{
	Comparison of the optimized inputs applied on each model. OT: Inputs generated with OC based on \textit{TOMGRO}, OS: Inputs generated with OC based on \textit{SIMPLE}. I.e., OT-SIMPLE is the biomass predicted by the SIMPLE given the optimized inputs generated by the \textit{TOMGRO}. 
	} 
	\label{fig:dssat_opt_comp}
	\vspace{-0.5cm}
\end{figure}
\label{sec:compare-optimal-results}

\vspace{-0.1cm}
\section{Discussion}
\label{sec:discussion}
Overall, all three models achieve similar yields under standard weather and \ac{GH} conditions, even though only one of the models, the reduced \textit{TOMGRO}, has initially been designed for \ac{GH} applications. The existing \textit{DSSAT} implementation rendered it inaccessible for immediate use in the current implementation of our \ac{OC} approaches. The remaining two models have different structures and states, with \textit{TOMGRO} being more intuitive to measure. All models share similar inputs, but \textit{DSSAT} allows for more degrees of freedom, while \textit{SIMPLE} allows for the least.

Based on the optimization results, the crop models have different sensitivities to environmental variables. We applied the \ac{OC} trajectories of the \ac{GH} to both models and found that while the resulting biomass depends on the chosen trajectory, the differences between the models are not substantial, even when the \ac{GH} is optimized with the other model, shown in Fig.~\ref{fig:dssat_opt_comp}. Comparing to \textit{DSSAT}, its prediction lies in between the other two models, leaving uncertainty as to which prediction is closer to reality. These findings emphasize the importance of thorough model calibration, as accurate validation does not necessarily lead to similar optimal \ac{GH} environment trajectories for different crop models.

\begin{table}[tb]
\vspace{-0.0cm}
\centering
\caption{Overall Comparison of the three tomato models. \\ \,}
\label{tab:Summary}
\scalebox{0.77}{
\begin{tabular}{l|c|c|c|c|c|c|c|c}

\cline{2-9}

& \multicolumn{1}{c|}{\rotatebox[origin=c]{90}{\textbf{\begin{tabular}[c]{@{}c@{}}Developed for \\ \acp{GH}\end{tabular}}}}
& \multicolumn{1}{c|}{\rotatebox[origin=c]{90}{\textbf{\begin{tabular}[c]{@{}c@{}}Valid under\\ \ac{GH} conditions\end{tabular}}}}
& \multicolumn{1}{c|}{\rotatebox[origin=c]{90}{\textbf{\begin{tabular}[c]{@{}c@{}}Usable for \\ \ac{OC} \end{tabular}}}}
& \multicolumn{1}{c|}{\rotatebox[origin=c]{90}{\textbf{\begin{tabular}[c]{@{}c@{}}Optimal \\ harvest in \si{\kilogram\per\metre\squared}\end{tabular}}}}   
& \multicolumn{1}{c|}{\rotatebox[origin=c]{90}{\textbf{\begin{tabular}[c]{@{}c@{}}Economic \\ yield in \si{\sieuro\per\square\metre} \end{tabular}}}}
& \multicolumn{1}{c|}{\rotatebox[origin=c]{90}{\textbf{\begin{tabular}[c]{@{}c@{}}Optimal harvest \\ based on \textit{DSSAT} \\in \si{\kilogram\per\metre\squared}\end{tabular}}}} 
& \multicolumn{1}{c|}{\rotatebox[origin=c]{90}{\textbf{\begin{tabular}[c]{@{}c@{}}Economic yield \\ based on \textit{DSSAT} \\ in \si{\sieuro\per\square\metre} \end{tabular}}}}
& \multicolumn{1}{c|}{\rotatebox[origin=c]{90}{\textbf{\begin{tabular}[c]{@{}c@{}} \, Computational \, \\ cost in \si{\second} \end{tabular}}}}   
                                                               \\ 

\hline

\multicolumn{1}{|l|}{\textbf{SIMPLE}}                                                            
& \multicolumn{1}{c|}{\redcross}                         
& \multicolumn{1}{c|}{\greencheck}                              
& \multicolumn{1}{c|}{\greencheck}                          
& \multicolumn{1}{c|}{16.73}   
& \multicolumn{1}{c|}{16.16}
& \multicolumn{1}{c|}{17.19}    
& \multicolumn{1}{c|}{16.78}
& \multicolumn{1}{c|}{\textbf{14}}                                                                           
                       \\ 

\hline 

\multicolumn{1}{|l|}{\textbf{\begin{tabular}[l]{@{}l@{}}Reduced\\TOMGRO \end{tabular}}}         
& \multicolumn{1}{c|}{\greencheck}   
& \multicolumn{1}{c|}{\greencheck } 
& \multicolumn{1}{c|}{\greencheck}                          
& \multicolumn{1}{c|}{19.76}     
& \multicolumn{1}{c|}{21.88}  
& \multicolumn{1}{c|}{17.43} 
& \multicolumn{1}{c|}{\textbf{18.02}}   
& \multicolumn{1}{c|}{19}                                       
 \\ 

\hline

\multicolumn{1}{|l|}{\textbf{DSSAT}}                                                            
& \multicolumn{1}{c|}{\redcross}                         
& \multicolumn{1}{c|}{\greencheck}                              
& \multicolumn{1}{c|}{\redcross }                          
& \multicolumn{1}{c|}{-}                                        
& \multicolumn{1}{c|}{-}                                       
& \multicolumn{1}{c|}{-}
& \multicolumn{1}{c|}{-}      
& \multicolumn{1}{c|}{-}                                                                           

                       \\ 

\hline 
\end{tabular}
}
\vspace{-0.5cm}
\end{table}

The results are presented in Table \ref{tab:Summary}, where it is shown that the \textit{SIMPLE} model has a computation time 1.35 times faster than the \textit{TOMGRO} model on a standard laptop, given an identical parallelized implementation in Python. 
This outcome holds noteworthy implications for closed-loop control systems that depend on updating constrained optimization using measured states. While the current \ac{GH} model only permits hourly updates, adapting to a finer granularity would be advantageous for resource optimization, making the differences in computation time more consequential. Furthermore, the current configuration does not account for spatially distributed control approaches that consider spatial discrepancies in the \ac{GH} and crop states, which would lead to a substantial increase in computational resource and real-time capability requirements.
In order to compensate for the explained model parameter uncertainties and to enable a fair comparison, we compare the achievable harvest and yield using the \textit{DSSAT} model predictions, using the inputs generated in the optimization with both \textit{SIMPLE} and \textit{TOMGRO} (see columns 6 and 7). Although both models achieve similar final biomass, the \textit{TOMGRO} model produces slightly better economic yield (7 \% higher). However, in a closed-loop regime, small model inaccuracies are less important due to constant feedback.  In summary, both models have their strengths and weaknesses, but the \textit{SIMPLE} model seems more attractive for closed-loop applications since it can be easily adapted to other crops. Nevertheless, obtaining the states poses difficulties that must be overcome to use the \textit{SIMPLE} model in a closed-loop.

\vspace{-0.1cm}
\section{Conclusion}
\label{sec:conclusion}
In this study, we compared three different models for the growth of a tomato plant, namely the \textit{SIMPLE}, reduced \textit{TOMGRO} and the \textit{DSSAT} model. All crop models can achieve similar accurate yields under standard weather and \ac{GH} conditions. While the \textit{DSSAT} model turned out to be inaccessible for optimal control approaches, the remaining two models, \textit{TOMGRO} and \textit{SIMPLE}, were compared in a combination with the same \ac{GH} model with respect to the usability in an \ac{OC} approach. The optimization results revealed that the models have different sensitivities to environmental variables, emphasizing the importance of careful model calibration. The \textit{SIMPLE} model had a faster computation time, making it more attractive for closed-loop control systems, but the \textit{TOMGRO} model achieved an overall higher yield, leaving uncertainty about which model is closer to reality. Despite their strengths and weaknesses, both models are relevant for closed-loop applications. The \textit{SIMPLE} model is particularly useful due to its adaptability to other crops. However, challenges associated with obtaining states need to be overcome for the use in a closed-loop control system which will be left for future work.
\vspace{-0.3cm}

\bibliography{Agreta2023}
\bibliographystyle{unsrt}

\end{document}